\documentclass [11pt,twoside,a4paper]{article}
\usepackage{amsfonts}
\usepackage{amsmath}
\usepackage{amstext}
\usepackage{amssymb}
\usepackage{mathrsfs}
\usepackage{amscd}
\usepackage{xypic}
\usepackage{enumitem}           
\usepackage{epsf}              
\usepackage{graphicx}          
\usepackage{fancybox}          
\usepackage{color}             
\usepackage{fancyhdr}
\usepackage[hang,footnotesize]{caption2}  

\def\mathcal{\mathscr}
\newfont{\aaa}{cmb10 at 19pt}
\newfont{\bbb}{cmb10 at 11pt}
\newtheorem{thm}{Theorem}[section]

\newtheorem{lem}[thm]{Lemma}
\newtheorem{prop}[thm]{Proposition}

\newtheorem{rem}[thm]{Remark}

\def\le{\leqslant}
\def\leq{\leqslant}
\def\ge{\geqslant}
\def\geq{\geqslant}

\def\l{{\langle}}
\def\r{\rangle}

\def\R{{\mathbb R}}
\def\C{{\mathbb C}}

\def\a{\alpha}

\def\eps{\varepsilon}

\def\E{{\mathbb E}}
\def\P{{\mathbb P}}
\def\Z{{\mathbb Z}}
\def\N{{\mathbb N}}

\newcommand{\beq}{\begin{equation}}
\newcommand{\eeq}{\end{equation}}
\newcommand{\bey}{\begin{eqnarray}}
\newcommand{\eey}{\end{eqnarray}}
\newcommand{\beyy}{\begin{eqnarray*}}
\newcommand{\eeyy}{\end{eqnarray*}}

\makeatletter
\def\@evenhead{
\vbox{\hbox to \textwidth {}{\hspace{0mm}{\footnotesize
\thepage}}{\hspace{8.5cm} {\footnotesize {Z. CHEN, D. WU, Y. XIAO}}} \protect\vspace{1truemm}\relax \hrule depth0pt
height0.15truemm width\textwidth}}
\def\@evenfoot{}
\def\@oddhead{\vbox{\hbox to \textwidth
{{\hspace{0cm}{\footnotesize Smoothness of the local times of anisotropic Gaussian random fields}\hfill{\footnotesize
\thepage}}\hspace{0mm}}{} \protect\vspace{1truemm}\relax\hrule
depth0pt height0.15truemm width\textwidth}}
\def\@oddfoot{}
\makeatother

\newcommand{\pf}{\noindent{\it Proof}\quad}

\begin{document}

\thispagestyle{empty} \thispagestyle{fancy} {
\fancyhead[lO,RE]{\footnotesize  Front. Math. China \\
DOI xxxxxx\\[3mm]
\fancyhead[RO,LE]{\scriptsize \bf 
} \fancyfoot[CE,CO]{}}
\renewcommand{\headrulewidth}{0pt}


\setcounter{page}{1}
\qquad\\[8mm]

\noindent{\aaa{Smoothness of Local Times and Self-Intersection Local Times of Gaussian Random Fields}}\\[1mm]

\noindent{\bbb Zhenlong CHEN$^1$,\quad Dongsheng WU$^2$, \quad Yimin XIAO$^3$}\\[-1mm]

\noindent\footnotesize{1 School of Statistics and Mathematics, Zhejiang Gongshang University, China}

\noindent\footnotesize{2 Department of Mathematical Sciences, University of Alabama in Huntsville, U.S.A.}

\noindent\footnotesize{3 Department of Statistics and Probability, Michigan State University, U.S.A.}\\[6mm]

\vskip-2mm \noindent{\footnotesize$\copyright$ Higher Education
Press and Springer-Verlag Berlin Heidelberg 2013} \vskip 4mm

\normalsize\noindent{\bbb Abstract}\quad
This paper is concerned with the smoothness (in the sense of Meyer-Watanabe) of the local times
of Gaussian random fields. Sufficient and necessary conditions for the existence and smoothness
of the local times, collision local times, and self-intersection local times are established
for a large class of Gaussian random fields, including fractional Brownian motions, fractional Brownian sheets
and solutions of stochastic heat equations driven by space-time Gaussian noise.
\vspace{0.3cm}

\footnotetext{Received xxxx; accepted xxxx\\
\hspace*{5.8mm}Corresponding author: Yimin XIAO, E-mail: xiao@stt.msu.edu}

\noindent{\bbb Keywords}\quad Anisotropic Gaussian field, local time,
collision local time, intersection local time, self-intersection local time,
chaos expansion\\
{\bbb MSC}\quad 60G15, 60H05, 60H07\\[0.4cm]

\setcounter{equation}{0}
\setcounter{section}{1}
\setcounter{thm}{0}

\noindent{\bbb 1\quad Introduction}\label{Sec:Intro}\\[0.1cm]

In recent years, Malliavin calculus has been shown to be very useful in stochastic
analysis of Gaussian processes [cf. Nualart (2006)]. In particular, many authors
have studied the chaos expansion and smoothness in the sense of Meyer-Watanabe
of local times and intersection local times of Brownian motion, fractional Brownian
motion and related self-similar Gaussian processes. See Nualart and Vives (1992),
Imkeller and Weisz (1995, 1999), Imkeller et al. (1995), Hu (2001), Hu and
\O ksendal (2002), Eddahbi and Vives (2003), Hu and Nualart (2005),  Yan et al.
(2009), Jiang and Wang (2009), Yan and Shen (2010),  Chen and Yan (2011), Shen
and Yan (2011), Shen el al. (2012). However,  there have been only a
few results on smoothness of local times of Gaussian random fields due to their more
complicated dependence structures. We refer to Imkeller and Weisz (1994, 1999) for the
case of Brownian sheet and to Eddahbi, et al. (2005, 2007) for  results on fractional 
Brownian sheets.

The main purpose of this paper is to study the smoothness in the sense of Meyer-Watanabe of
the local times of a large class of Gaussian random fields, including fractional Brownian sheets and
solutions of stochastic heat equations driven by space-time Gaussian noise.
More specifically,   let $X=\{X(t), t \in \R^N\}$ be a Gaussian random field with
values in $\R^d$ defined on a probability space $(\Omega, {\cal F}, \P)$ by
\begin{equation}\label{def:X2}
X(t) = \big(X_1(t), \ldots, X_d(t)\big), \quad \forall\, t \in
\R^N.
\end{equation}
We will call $X$ an $(N, d)$-Gaussian random field. We assume that
the coordinate fields $X_1, \ldots, X_d$ are independent copies
of a real-valued, centered Gaussian random field $X_0 = \{X_0(t),
t \in \R^N\}$ with continuous covariance function $R(s, t) = \E[X_0(s)X_0(t)]$.

Let $H= (H_1, \ldots, H_N) \in (0, 1)^N$ be a fixed vector. For $a, \, b
\in \R^N$ with $a_j < b_j$ ($j = 1, \ldots, N$), let $I = [a,\,b]:=\prod_{j=1}^N
[a_j, b_j] \subseteq \R^N$ be the compact interval (or a rectangle). For
simplicity, we will take $I = [0, 1]^N$ throughout this paper.
We further assume that $X_0 = \{X_0(t), t \in \R^N\}$ satisfies the following conditions:
\begin{itemize}
\item[(C1)]\ There exists a positive and finite constant $c_{1}$ such that
\begin{equation}\label{Eq:Var1}
\E\Big[\big(X_0(s) - X_0(t)\big)^2\Big] \le c_{1}\,
\sum_{j=1}^N |s_j - t_j|^{2 H_j}, \quad \forall \, s,\, t \in I.
\end{equation}
\item[(C2)]\ There exists a constant $c_{2} > 0$ such that
for all $s, t \in I$,
\begin{equation}\label{Eq:C2}
{\rm Var}\big(X_0(t)\big| X_0(s)\big) \ge c_{2}\,
\sum_{j=1}^N \min\Big\{|s_j - t_j|^{2 H_j}, |t_{j}|^{2 H_j}\Big\}.
\end{equation}
Here ${\rm Var}(X_0(t)|X_0(s))$ denotes the conditional variance
of $X_0(t)$ given $X_0(s)$.
\end{itemize}

The class of Gaussian random fields that satisfy Conditions (C1) and (C2) is large.
When $N=1$, it includes fractional Brownian motion, bi-fractional Brownian motion
and related Gaussian processes. For $N \ge 2$, this class contains  fractional
Brownian sheets [cf. Ayache and Xiao (2005), Wu and Xiao (2007) for verification],
solutions to stochastic heat equation driven by space-time Gaussian noises
[Mueller and Tribe (2002), Wu and Xiao (2006), Dalang, et al.  (2015), Tudor and Xiao (2015)] and
many more [cf. Xiao (2009)].

The purpose of this paper is to study the existence and  smoothness (in the
sense of Meyer-Watanabe) of the local times and the self-intersection local
times of Gaussian random fields that satisfy Conditions (C1), (C2) and/or
(C3) below. Our
main results in Sections 2 and 3 unify and extend the previous results in
the references mentioned at the beginning of the Introduction. We should also
mention that H\"older regularities of local times and their applications
to sample path properties of Gaussian random fields have been studied by
several authors, including Pitt (1978), Geman and Horowitz (1980), Xiao
and Zhang (2002), Ayache and Xiao (2005), Ayache, et al. (2008), Wu and
Xiao (2007, 2010, 2011), Xiao (2009), Bierm\'{e} et al. (2009), Chen and
Xiao (2012).

The rest of this paper is organized as follows. In Section 2 we provide
a sufficient and necessary condition for the existence, and a sufficient
condition for the smoothness (in the sense of Meyer-Watanabe) of the
local time at any level $x \in \R^d$ for a large class of Gaussian random
fields. We also prove that this condition for the smoothness is also
necessary for the local times at $x= 0$. We then apply the conditions to prove the
existence and smoothness results for the collision local times and
the intersection local times for two independent anisotropic
Gaussian random fields.

Section 3 is concerned with self-intersection local times. We establish
a sufficient and necessary condition for the existence and smoothness of
self-intersection local times on two disjoint intervals. More interestingly,
we also consider the analogous problems on two intersecting intervals.
We will see that the results in the intersecting cases are different from
and more difficult than those in the disjoint case.


Throughout this paper, we will use $c$ to denote unspecified
positive finite constants which may be different in each
appearance.  More specific constants  are numbered
as $c_{1},\ c_{2},\ldots$.

\noindent\\[1mm]

\setcounter{equation}{0}
\setcounter{section}{2}
\setcounter{thm}{0}
\noindent{\bbb 2\quad Existence and smoothness of the local times}\\[0.1cm]

This section is concerned with the existence and smoothness of the
local times of a Gaussian random field $X$ in the sense of Meyer-Watanabe.
We start by recalling the definition of Chaos expansion, which is an orthogonal
decomposition of $L^{2}(\Omega, \mathbb{P})$. We refer to Nualart and Vives
(1992), Meyer (1993), Hu (2001), Nualart (2006) and the references therein
for more information.

Let $\Omega$ be the space of continuous $\R^d$-valued functions
$\omega$ on $I.$ Then $\Omega$ is a Banach space with respect to
the sup norm. Let $\mathcal{F}$ be the Borel $\sigma$-algebra on
$\Omega.$ Let $\P$ be a probability measure on
$(\Omega,\mathcal{F})$, and $\E$ denote the expectation on this
probability space. Denote by $L^2(\Omega, \mathbb{P})$ the space of all
 real (or complex) valued functional on
$\Omega$ such that
\[\E(F^2)=\int_\Omega|F(\omega)|^2\P(d\omega)<\infty.\]

Let  $Y=\left\{Y_1(t),\ldots,Y_d(t),\,t\in I\right\}$ be an
$(N,d)$-Gaussian random field, where $Y_1,\ldots, Y_d$ are $d$
independent copies of some centered, real-valued Gaussian random
field $Y_0$ on $I$.
Let $p_n(y_1,\ldots,y_k)$ be a polynomial of degree $n$ of $k$
variables $y_1,\ldots,y_k$. Then, for any $t^1,\ldots,t^k\in I$ and
$i_1,\ldots,i_k\in\{1,\ldots,d\}$, $p_n\big(Y_{i_1}(t^1),\ldots,
Y_{i_k}(t^k)\big)$ is called a {\it
polynomial functional} of $Y$. Let ${\mathcal P}_n$ be the
completion with respect to the $L^2(\Omega, \mathbb{P})$ norm of
the set of all polynomials of degree less than or equal to $n$.  Then
${\mathcal P}_n$ is a subspace of $L^2(\Omega, \mathbb{P})$.
Let ${\mathcal C}_n$ be the orthogonal
complement of ${\mathcal P}_{n-1}$ in ${\mathcal P}_n.$ Then
 $L^2(\Omega, \mathbb{P})$ is the direct sum of ${\mathcal C}_n,$ i.e.,
 \[L^2(\Omega, \mathbb{P})=\bigoplus_{n=0}^\infty {\mathcal C}_n.\]
Namely, for any functional $F\in L^2(\Omega, \mathbb{P})$, there
exists a sequence $\{F_n\}_{n=0}^\infty$ with $F_n\in{\mathcal C}_n,$  such
that $F=\sum_{n=0}^\infty F_n.$
This decomposition is called the {\it chaos expansion} of $F$, and
$F_n$ is called the {\it $n$-th chaos} of $F$.  Clearly,
\[F_0=\E (F), \quad \E(|F|^2)=\sum_{n=0}^\infty\E(|F_n|^2).\]

In Malliavin Calculus, the space of ``smooth'' functions in
the sense of Meyer-Watanabe [cf. Watanabe (1984), Nualart (2006)] is defined by
\[
D_1:=\bigg\{F\in L^2(\Omega, \mathbb{P}),\,F=\sum_{n=0}^\infty F_n \ \mbox{ and }\
\sum_{n=0}^\infty n\E(|F_n|^2)<\infty\bigg\}.
\]
For $F\in L^2(\Omega, \mathbb{P})$ with a chaos expansion $F=\sum
F_n$, define the operator $\Gamma_u$ with $u\in[0,\,1]$ by
\begin{equation}\label{Def:GammauF}
\Gamma_u F:=\sum_{n=0}^\infty u^nF_n,
\end{equation}
 and set
$\Theta_F(u):=\Gamma_{\sqrt u}F.$ Clearly $\Theta_F(1)=F$.  Define
$\Phi_{\Theta_F}(u):=\frac{d}{du}\E(|\Theta_F(u)|^2),$
we then have
\[\Phi_{\Theta_F}(u)=\sum_{n=1}^\infty n u^{n-1}\E(|F_n|^2).\]

In the following, we provide several technical lemmas which will be useful
 for proving the existence and smoothness of local times.
Lemma \ref{lem:lemma1} is similar to Lemma 8.6 in Bierm\'e, Lacaux
and Xiao (2009) whose proof is elementary.  Lemmas
\ref{lem:lemma2} and \ref{lem:lemma3} are  from Wu and Xiao (2010).

\begin{lem}\label{lem:lemma1}
 Let $\alpha$ and $\beta$  be positive constants, then for all $A\in (0,1)$
\begin{equation}\label{Eq:lemma1}
\arraycolsep=1pt\begin{array}{ll} \displaystyle
\int_{0}^{1}\frac{1}{\big(A+t^{\alpha}\big)^{\beta} }dt\asymp
\left\{\begin{array}{ll}
A^{-(\beta-\frac{1}{\a}}) \ \ \ \qquad &\hbox{ if } \, \alpha\beta > 1,\\
\log \big(1+A^{-\frac{1}{\a}} \big) &\hbox{ if }\, \alpha\beta = 1,\\
1 &\hbox{ if }\, \a \beta< 1.
\end{array}
\right.
\end{array}
\end{equation}
  In the above, $f(A)\asymp g(A)$  means that the ratio $ f (A)/g(A)$ is bounded from
below and above by positive constants that do not depend on $A\in
(0, 1)$.
\end{lem}

\begin{lem}\label{lem:lemma2}
Let $\alpha$ and $\beta$  be positive constants such that $\alpha\beta\geq
 1$.
\begin{itemize}
\item[\rm (i)]\  If $\alpha\beta> 1$, then there exists a constant $c_3> 0$  whose
value depends on $\alpha$ and $\beta$ only such that for all $A\in (0, 1),\ r >
 0,\  u^{*}\in\mathbb{R}$,   all
integers $n\geq 1$ and all distinct $u_{1},\ldots, u_{n}\in
O(u^{*}, r)$ we have
\begin{equation}\label{Eq:lemma21}
\arraycolsep=1pt\begin{array}{ll} \displaystyle\int_{O(u^{*},
r)}\displaystyle\frac{du}{\big(A+\displaystyle\min_{1\leq j\leq
n}|u-u_{j}|^{\alpha}\big)^{\beta} }\leq c_3\,
nA^{-(\beta-\frac{1}{\a}}).
\end{array}
\end{equation}
where $O(u^{*}, r)$ denotes a ball centered at $u^{*}$ with radius
$r$.
\item[\rm (ii)]\ If $\alpha\beta=1$, then for any $\kappa\in (0,1)$ there
exists a constant $c_4> 0$   whose value depends on $\alpha$, $\beta$
and $\kappa$ only such that for all $A\in (0, 1),\ r > 0,\
 u^{*}\in\mathbb{R}$,   all
integers $n\geq 1$ and all distinct $u_{1},\ldots, u_{n}\in
O(u^{*}, r)$ we have
\begin{equation}\label{Eq:lemma22}
\arraycolsep=1pt\begin{array}{ll} \displaystyle \int_{O(u^{*},
r)}\frac{du}{\big(A+\displaystyle\min_{1\leq j\leq
n}|u-u_{j}|^{\alpha}\big)^{\beta} }\leq
c_4\, n\log\bigg[e+\bigg(\frac{r}{n}A^{-\frac{1}{\a}}\bigg)^{\kappa}\bigg].
\end{array}
\end{equation}
\end{itemize}
\end{lem}

\begin{lem}\label{lem:lemma3}
Let $0<\beta<1$  be a constant. Then there exists a positive constant $c_5$ such
that the following statements hold.
\begin{itemize}
\item[\rm (i)]\  For all $ r > 0,\
 u^{*}\in\mathbb{R}$,   all
integers $n\geq 1$ and all distinct $u_{1},\ldots, u_{n}\in
O(u^{*}, r)$ we have
\begin{equation}\label{Eq:lemma31}
\arraycolsep=1pt\begin{array}{ll} \displaystyle \int_{O(u^{*},
r)}\frac{du}{\displaystyle\min_{1\leq j\leq n}|u-u_{j}|^{\beta}
}\leq c_5\, n^{\beta}r^{-(\beta-1)}.
\end{array}
\end{equation}
\item[\rm (ii)]\ For all constants
$r > 0$ and $M>0$, all $u^{*}\in\mathbb{R}$,  integers $n\geq 1$
and all distinct $u_{1},\ldots, u_{n}\in O(u^{*}, r)$ we have
\begin{equation}\label{Eq:lemma32}
 \int_{O(u^{*}, r)}\log\Big[e+M\big(\min_{1\leq j\leq
n}|u-u_{j}|\big)^{-\beta}\Big]\, du\leq
c_5\, r\log\bigg[e+M\Big(\frac{r}{n}\bigg)^{-\beta}\Big].
\end{equation}
\end{itemize}
\end{lem}

\vskip 0.1cm

\noindent{\bbb 2.1\quad General results }\\[0.1cm]

We will apply the following proposition and the method of its
proof to study the existence and smoothness of the local times of $X$.

\begin{prop}\label{prn:prn}
Let $X=\{X(t),\,t\in I\}$ be an $(N,d)$-Gaussian field
defined by (\ref{def:X2}) and assume that $X_{0}$ satisfies
Conditions (C1) and (C2) with index $H\in (0, 1)^N$.  Then, for any
$\gamma>0$, $\lambda\geq 0$,
\begin{equation}\label{Eq:prn}
\arraycolsep=1pt\begin{array}{ll} \displaystyle
\int_{I^{2}}\frac{\left|\big[\mathbb{E}\big(X_0(s)X_0(t)\big)\big]\right|^{\lambda}}{[\det{\rm
Cov}(X_{0}(s),  X_{0}(t))]^{\frac{\gamma}{2}} }\, dsdt<\infty
\end{array}
\end{equation}
if and only if $\sum_{\ell=1}^{N}1/H_{\ell}>\gamma$.
\end{prop}

\pf First we prove the sufficiency. By (C2) we have
\begin{equation}\label{Eq:prn0}
{\rm Var}(X_0(s))\geq {\rm  Var}\left(X_0(s)|X_0\big(\frac{s}{2}\big)\right)\ge
c_{2}2^{-2}\sum_{j=1}^N s_j^{2H_j}, \quad  \forall \, s \in I.
\end{equation}
This and  the fact that
\begin{equation}\label{Eq:CovVar}
{\rm detCov}\left(X_0(s), X_0(t)\right)={\rm Var}(X_0(s)){\rm
Var}\left(X_0(t)|X_0(s)\right)
\end{equation}
imply
\begin{equation}\label{Eq:prn1}
\begin{split}
&{\rm detCov}\left(X_0(s), X_0(t)\right)\ge c\, \bigg(\sum_{j=1}^N
s_j^{2H_j}\bigg) \bigg(\sum_{j=1}^N \min\Big\{|s_j - t_j|^{2 H_j},
t_{j}^{2 H_j}\Big\}\bigg).
\end{split}
\end{equation}
On the other hand, it follows from the Cauchy-Schwarz inequality
and the continuity of the covariance function $R(s, t)$ that
 \begin{equation}\label{Eq:EXX}
 \left|\big[\mathbb{E}\big(X_0(s)X_0(t)\big)\big]\right|^{\lambda}\leq
 c, \quad \forall\ s,t\in I.
\end{equation}
Hence, for proving the sufficiency, it suffices to verify that
if $\sum_{j=1}^{N}\frac1{H_j}>\gamma$, then
\begin{equation}\label{Eq:prn2}
\int_{I^2}\frac{ds\,dt }{\big[\sum_{j=1}^Ns_j^{2H_j}\big]^{\frac \gamma
2}\big[\sum_{j=1}^N \min\big\{|s_j - t_j|^{2 H_j}, t_{j}^{2
H_j}\big\}\big]^{\frac \gamma 2}} <\infty.
\end{equation}

To estimate the  integral in  (\ref{Eq:prn2}), we will assume that
\begin{equation}\label{Eq:H}
0< H_1\leq  H_2\leq\cdots\leq  H_N<1
\end{equation}
and integrate in the order of $dt_{1}, \ldots, dt_{N}, ds_{1}, \ldots,
ds_{N}$. When $\sum_{j=1}^N\frac1{H_j}>\gamma$, there exists
$k\in \{1, 2, \ldots, N\}$ such that
\begin{equation}\label{Eq:prn3}
\sum_{j=1}^{k-1}\frac1{H_j}\leq \gamma<\sum_{j=1}^{k}\frac1{H_j}.
\end{equation}
Note that
\begin{equation}\label{Eq:prn4}
\begin{split}
&\int_{I}\frac{dt_{1} \ldots dt_{N}}{\big[\sum_{j=1}^N
\min\big\{|s_j - t_j|^{2 H_j}, t_{j}^{2 H_j}\big\}\big]^{\frac
\gamma 2}} \\
& \qquad \leq \int_{I}\frac{dt_{1} \ldots dt_{N}}{\big[\sum_{j=1}^k
\min\big\{|s_j - t_j|^{2 H_j}, t_{j}^{2 H_j}\big\}\big]^{\frac
\gamma 2}}.
 \end{split}
\end{equation}
We distinguish two cases: (i) $\sum_{j=1}^{k-1}\frac1{H_j}<
\gamma<\sum_{j=1}^{k}\frac1{H_j}$ and (ii)
$\sum_{j=1}^{k-1}\frac1{H_j}= \gamma<\sum_{j=1}^{k}\frac1{H_j}$,
and show that the last integral in
(\ref{Eq:prn4}) is bounded by a constant that is independent of $s \in I$.

In Case (i), if $k=1$, then  $H_{1}\gamma<1$. We apply (i) of
Lemma \ref{lem:lemma3} to derive
\begin{equation*}\label{Eq:prn5}
\begin{split}
\int_{I}\frac{dt_{1} \ldots dt_{N}}{\big[\sum_{j=1}^k
\min\big\{|s_j - t_j|^{2 H_j}, t_{j}^{2 H_j}\big\}\big]^{\frac
\gamma 2}}
\leq\int_{I}\frac{dt_{1} \ldots dt_{N}}{\big[ \min\big\{|s_1 -
t_1|^{2 H_1}, t_{1}^{2 H_1}\big\}\big]^{\frac \gamma 2}} \le c_6,
\end{split}
\end{equation*}
as desired. If $k>1$, then $H_{1}\gamma>1$. We first apply (i) of Lemma
\ref{lem:lemma2}
 with $\alpha=2H_{1}$, $\beta=\frac{\gamma}{2}$  and $A= \sum_{j=2}^N
\min\big\{|s_j - t_j|^{2 H_j}, t_{j}^{2 H_j}\big\}$ to deduce that
\begin{equation}\label{Eq:prn6}
\begin{split}
&\int_{0}^{1}\frac{dt_{1} }{\Big[\min\big\{|s_1 - t_1|^{2 H_1},
t_{1}^{2 H_1}\big\}+ \sum_{j=2}^k \min\big\{|s_j - t_j|^{2 H_j},
t_{j}^{2 H_j}\big\}\Big]^{\frac
\gamma 2}}\\
&\qquad \leq \frac{c_7}{\big[\sum_{j=2}^k \min\big\{|s_j - t_j|^{
2H_j}, t_{j}^{
2H_j}\big\}\big]^{\frac{1}{2}(\gamma-\frac{1}{H_{1}})}},
 \end{split}
\end{equation}
where $c_7$ is a constant which only depends on $H_{1}$ and $\gamma$. By
repeatedly using Part (i) of Lemma \ref{lem:lemma2} as in
(\ref{Eq:prn6}), after $k-1$ steps, we obtain that
\begin{equation}\label{Eq:prn7}
\begin{split}
&\int_{I}\frac{dt_{1} \ldots dt_{N}}{\big[\sum_{j=1}^k
\min\big\{|s_j - t_j|^{2 H_j}, t_{j}^{2 H_j}\big\}\big]^{\frac
\gamma 2}}\\
& \qquad \leq c\, \int_{0}^{1}\frac{ dt_{k} }{\big[
\min\big\{|s_k - t_k|^{ 2H_k}, t_{k}^{
2H_k}\big\}\big]^{\frac{1}{2}(\gamma-\sum_{j=1}^{k-1}\frac{1}{H_{j}})}}.
 \end{split}
\end{equation}
Notice that $H_{k}(\gamma-\sum_{j=1}^{k-1}\frac{1}{H_{j}})<1$, by
applying (i) of Lemma \ref{lem:lemma3} to the last integral in
(\ref{Eq:prn7}), we see from (\ref{Eq:prn4}) that in Case (i)
\begin{equation}\label{Eq:prn8}
\begin{split}
\int_{I}\frac{dt_{1} \ldots dt_{N}}{\big[\sum_{j=1}^N
\min\big\{|s_j - t_j|^{2 H_j}, t_{j}^{2 H_j}\big\}\big]^{\frac
\gamma 2}} \le c_8.
 \end{split}
\end{equation}

Now we consider Case (ii). Notice that $k>1$ in (\ref{Eq:prn4}).  We integrate in order
of $dt_{1},  \ldots, dt_{N}$  and  repeatedly apply Part (i) of Lemma
\ref{lem:lemma2}   for  $k-2$ steps  to get
\begin{equation}\label{Eq:pron1}
\begin{split}
&\int_{I}\frac{dt_{1} \ldots dt_{N}}{\big[\sum_{j=1}^k
\min\big\{|s_j - t_j|^{2 H_j}, t_{j}^{2 H_j}\big\}\big]^{\frac
\gamma 2}}\\
&\leq c\, \int_{0}^1 \int_0^1\frac{dt_{k-1} \, dt_{k}
}{\big[\sum_{j=k-1}^{k} \min\big\{|s_j - t_j|^{ 2H_j},
t_{j}^{
2H_j}\big\}\big]^{\frac{1}{2}(\gamma-\sum_{j=1}^{k-2}\frac{1}{H_{j}})}}.
 \end{split}
\end{equation}
Note that $H_{k-1} (\gamma-\sum_{j=1}^{k-2}\frac{1}{H_{j}})=1$. By applying  (ii)
of Lemma \ref{lem:lemma2} with $A=\min\big\{|s_k - t_k|^{ 2H_k},
t_{k}^{ 2H_k}\big\}$ and Part (ii) of Lemma \ref{lem:lemma3}, we derive
\begin{equation}\label{Eq:pron2}
\begin{split}
&\int_{0}^{1}\int_{0}^{1}\frac{ dt_{k-1} dt_{k}
}{\big[\sum_{j=k-1}^{k} \min\big\{|s_j - t_j|^{2 H_j}, t_{j}^{
2H_j}\big\}\big]^{\frac{1}{2}(\gamma-\sum_{j=1}^{k-2}\frac{1}{H_{j}})}}\\
&\leq
c\, \int_{0}^{1}\log\Bigg[e+\Bigg(\frac{1}{2}\Big(\min\big\{|s_k
- t_k|^{ 2H_k}, t_{k}^{
2H_k}\big\}\Big)^{-\frac{1}{2H_{k-1}}}\Bigg)^{\kappa}\Bigg]dt_{k}\\
&\leq c\,\log\Big[e+2^{H_{k}-H_{k-1}}\Big],
 \end{split}
\end{equation}
where $\kappa\in (0,1) $ is a constant and we have used the 
fact that $H_{k}\geq H_{k-1}$. It follows from  (\ref{Eq:pron1}) and (\ref{Eq:pron2})
that  (\ref{Eq:prn8}) also holds in Case (ii).

Hence, by (\ref{Eq:prn2}) and (\ref{Eq:prn8}), we have
\begin{equation*}\label{Eq:pron3}
\begin{split}
\int_{I^2}\frac{ ds\,dt}{\big[\sum_{j=1}^Ns_j^{2H_j}\big]^{\frac \gamma
2}\big[\sum_{j=1}^N \min\big\{|s_j - t_j|^{2 H_j}, t_{j}^{2
H_j}\big\}\big]^{\frac \gamma 2}}
\le c\, \int_{I}\frac{ds}{\big[\sum_{j=1}^Ns_j^{H_j}\big]^{
\gamma }}.
\end{split}
\end{equation*}
It is elementary to verify, by using Lemma \ref{lem:lemma1}, that
the last integral is finite provided $\sum_{j=1}^N\frac1{H_j}>\gamma$.
This proves (\ref{Eq:prn2}), and thus the sufficiency.

To prove the necessity, we prove that if
$\sum_{j=1}^N\frac1{H_j}\leq \gamma$ then
\begin{equation}\label{Eq:suf1}
\int_{I^2}\frac{\left|\big[\E\left(X_0(s)X_0(t)\right)\big]\right|^\lambda}{\big[{\rm
detCov}\left(X_0(s), X_0(t)\right)\big]^{\frac \gamma 2}}\,ds\,dt = \infty.
\end{equation}

For $\varepsilon_{0}\in(0,\,\frac{1}{2})$, let $I_{\varepsilon_{0}}:=[\varepsilon_{0},\,1]^N$.
Eq. (\ref{Eq:prn0}) and the uniform continuity of $R(s,t)$ on $I_{\varepsilon_{0}}^2$ imply
that there exists a constant $\delta_0>0$ such that for all $s, t\in [\varepsilon_{0},
\varepsilon_{0}+\delta_{0}]^{N}$,
\begin{equation}\label{Eq:suf3}
\arraycolsep=1pt\begin{array}{ll}
\mathbb{E}\big(X_{0}(s)X_{0}(t)\big) \geq
\frac{1}{2}\mathbb{E}\big(X^{2}_{0}(t)\big)\geq  c_9 > 0.
\end{array}
\end{equation}
On the other hand, it follows from (\ref{Eq:CovVar}) and  Condition (C1) that for all $s, t\in I$,
\begin{equation}\label{Eq:suf2}
\begin{split}
{\rm detCov}\big(X_0(s), X_0(t)\big)\le
c\, \sum_{j=1}^N|s_j-t_j|^{2H_j}.
\end{split}
\end{equation}
By (\ref{Eq:suf3}) and (\ref{Eq:suf2}), we derive
\begin{equation*}\label{Eq:suf5}
\begin{split}
\int_{I^2}\frac{\big|\E\left(X_0(s)X_0(t)\right)\big|^{\lambda}}
{\big[{\rm detCov}\left(X_0(s), X_0(t)\right)\big]^{\frac \gamma 2}}\,ds\,dt
\ge c\,\int_{[\varepsilon_{0},
\varepsilon_{0}+\delta_{0}]^{2N}}\frac{ds\,dt }{\big[\sum_{j=1}^N|s_j-t_j|^{H_j}\big]^{\gamma }}.
\end{split}
\end{equation*}
By using Lemma \ref{lem:lemma1} again, it is elementary to verify
that the last integral is infinite  when
$\sum_{j=1}^N\frac1{H_j}\leq \gamma$. This proves the necessity of
the proposition. \hfill{$\square$}

\vskip 0.2cm

In the following, we consider the existence of the local time of 
a Gaussian random field satisfying (C1) and (C2). Instead of using 
a Fourier analytic argument as in Xiao (2009) [see Geman and 
Horowitz (1980) for a systematic account], we approximate the Dirac 
delta function by the heat kernel
\begin{equation}\label{Def:pvx}
p_{\varepsilon}(x)=\frac{1}{(2\pi\varepsilon)^{d/2}}\exp\bigg(-\frac{\|x\|^{2}}
{2\varepsilon}\bigg),\quad  x\in\mathbb{R}^{d},
\end{equation}
and let
\begin{equation}\label{Def:selfv}
\arraycolsep=1pt\begin{array}{ll}
L_{\varepsilon}(x, I, X)&\displaystyle=\int_{I}p_{\varepsilon}(X(s)-x)ds\\
&=\displaystyle\frac{1}{(2\pi)^{d}}\int_{I}\int_{\mathbb{R}^{d}}\exp\Big(i\langle\xi,
X(s)-x\rangle-\frac{\varepsilon \|\xi\|^{2} }{2}\Big)d\xi ds.
\end{array}
\end{equation}

The following is a general result on existence of local times.

\begin{lem}\label{Lem:EY2}
Let $Y=\{\left(Y_1(t),\ldots,Y_d(t)\right),\,t\in I\}$ be an 
$(N,d)$-Gaussian random field, where $Y_1,\ldots, Y_d$ are $d$ 
independent copies of a centered, real-valued Gaussian random 
field $Y_0$ on $I$.  Then for any $y \in \R^d$, as
$\varepsilon\rightarrow 0$,  $L_{\varepsilon}(y, I, Y)$ 
converges to a limit $L(y, I, Y))$ in
$L^{2}(\Omega, \mathbb{P})$  if and only if
\begin{equation}\label{Eq:EY1}
\arraycolsep=1pt\begin{array}{ll}
\displaystyle\int_{I^{2}} \exp\bigg(-\frac{\|y\|^{2}\mathbb{E}[(Y_{0}(s)-Y_{0}(t))^{2}]}{{\rm
detCov}\left(Y_0(t), Y_0(s)\right)}\bigg) \frac{ dsdt}{\big[{\rm detCov}\left(Y_0(t),
Y_0(s)\right)\big]^{\frac d2}
}<\infty.
\end{array}
\end{equation}
\end{lem}
\pf\, Let $I_{2d}$ be the identity matrix of
order $2d$ and let $\Gamma_{\varepsilon, d }(s,t)=\varepsilon
I_{2d}+{\rm Cov}(Y(s),  Y(t))$.  For any $y \in \R^d$ and  $\varepsilon>0$,
Fubini's theorem and (\ref{Def:selfv}) imply
\begin{equation}\label{Eq:EY2}
\begin{split}
&\mathbb{E}(|L_{\varepsilon}(y, I, Y)|^{2})=\frac{1}{(2\pi)^{2d}}\int_{I^{2}} dsdt
\int_{\mathbb{R}^{2d}} e^{-\frac{\varepsilon}{2}( \|\xi\|^{2} + \|\eta\|^{2})}\\
& \qquad\qquad \qquad \times \mathbb{E}\exp\Big(i\langle\xi, Y(s)-y\rangle 
- i\langle\eta, Y(t)-y\rangle\Big)d\xi d\eta\\
&\displaystyle=\frac{1}{(2\pi)^{2d}}\int_{I^{2}}dsdt\int_{\mathbb{R}^{2d}}e^{-i\langle\xi -\eta,
y\rangle})  \exp\Big(-\frac{1}{2}(\xi,
\eta)\Gamma_{\varepsilon,d}(s,t)(\xi, \eta)^{T}\Big)d\xi d\eta\\
& =\frac{1}{(2\pi)^{2d}}\int_{I^{2}}
\exp\Big(-\frac{1}{2}(y, y)\Gamma^{-1}_{\varepsilon,d}(s,t)(y,
y)^{T}\Big) \frac{dsdt}{\sqrt{\det\Gamma_{\varepsilon,d}(s,t)}}.
\end{split}
\end{equation}
Since the coordinate processes of $Y$ are independent copies of
$Y_{0}$,  we have
\begin{equation}\label{Eq:EY3}
\arraycolsep=1pt\begin{array}{ll}
\displaystyle\det\Gamma_{\varepsilon,d}(s,t)&\displaystyle=
\big[\det\Gamma_{\varepsilon, 1 }(s,t)\big]^{d}
\end{array}
\end{equation}
and
\begin{equation}\label{Eq:Ey4}
\arraycolsep=1pt\begin{array}{ll}\displaystyle\frac{1}{2}(y,
y)\Gamma^{-1}_{\varepsilon,d}(s,t))(y,
y)^{T}=\frac{\|y\|^{2}
\big(2\varepsilon+\mathbb{E}[(Y_{0}(s)-Y_{0}(t))^{2}]\big)  } {\det\Gamma_{\varepsilon,1}(s,t)},
\end{array}
\end{equation}
where $\Gamma_{\varepsilon, 1 }(s,t)=\varepsilon I_{2}+{\rm
Cov}(Y_{0}(s), Y_{0}(t))$.  It follows from (\ref{Eq:EY3}), (\ref{Eq:Ey4}) and the
dominated convergence theorem that
\begin{equation}\label{Eq:EY5}
\begin{split}
 & \lim_{\varepsilon\rightarrow 0}\mathbb{E}(|L_{\varepsilon}(y, I, Y)|^{2})\\
&=\frac{1}{(2\pi)^{2d}}\int_{I^{2}}\exp\Big(-\frac{\|y\|^{2}
\mathbb{E}[(Y_{0}(s)-Y_{0}(t))^{2}]}{\det\Gamma_{0,1}(s,t)}\Big)\,
\frac{dsdt}{[\det\Gamma_{0,1}(s,
t )]^{d/2} }.
\end{split}
\end{equation}

Next we show that $\{L_{\varepsilon}(y, I, Y), \varepsilon>0\}$
is a Cauchy sequence in $ L^{2}(\Omega,
\mathbb{P})$ if and only if (\ref{Eq:EY1}) holds. For
all integers  $m,n\geq 1$, we assume, without loss of generality,
that $m=n+p$ for some integer $p$. Let $\Gamma_{n+p}(s,t) = (n+p
)^{-1}I_{2d} + {\rm Cov}\big(Y(s), Y(t)\big)$, $\Gamma_{n}(s,t) =
n^{-1}I_{2d} + {\rm Cov}\big(Y(s), Y(t)\big)$ and
$$
\Gamma_{ m+p,n}(s,t) = \bigg(\begin{array}{cc}
(m+p)^{-1}I_{d}    \quad & 0\\
0 & n^{-1}I_{d}
\end{array}
\bigg) + {\rm Cov}\big(Y(s), Y(t)\big).
$$
Then, it follows from Fubini's theorem and (\ref{Def:selfv}) that
\begin{equation*}\label{Eq:ECY1}
\begin{split}
&\mathbb{E}\Big[\big(L_{\frac{1}{n+p}}(y, I, Y)-L_{\frac{1}{n}}(y,
I, Y)\big)^{2}\Big] =\frac{1}{(2\pi)^{2d}}\int_{I^{2}}dsdt\int_{\mathbb{R}^{2d}}
e^{-i \langle\xi -\eta, y\rangle }\\
&\quad \times \bigg\{\exp\Big(-\frac{1}{2}(\xi,\eta)\Gamma_{n+p}(s,t)(\xi,\eta)^{T}\Big) +\exp\Big(-\frac{1}{2}(\xi,
\eta)\Gamma_{n}(s,t)(\xi,\eta)^{T}\Big)\\
&\qquad \qquad \quad -2\exp\Big(-\frac{1}{2}(\xi,\eta)\Gamma_{n+p,n}(s,t)(\xi,\eta)^{T}\Big)\bigg\}\xi d\eta\\
& =\frac{1}{(2\pi)^{2d}}\int_{I^{2}}\bigg\{\frac{1}{\sqrt{\det\Gamma_{n+p}(s,t)}}
\exp\Big(-\frac{1}{2}(y, y)\Gamma^{-1}_{n+p}(s,t)(y,
y)^{T}\Big)\\
&  \qquad \quad +\frac{1}{\sqrt{\det\Gamma_{n}(s,t)}}
\exp\Big(-\frac{1}{2}(y, y)\Gamma^{-1}_{n}(s,t)(y,
y)^{T}\Big)\\
& \qquad \quad -\frac{2}{\sqrt{\det\Gamma_{n+p,,n}(s,t)}}
\exp\Big(-\frac{1}{2}(y, y)\Gamma^{-1}_{n+p,n}(s,t)(y,
y)^{T}\Big)\bigg\}dsdt.
\end{split}
\end{equation*}
Similarly to  (\ref{Eq:EY5}), we can verify that
\begin{equation*}\label{Eq:part21}
\arraycolsep=1pt\begin{array}{ll} &\displaystyle\lim_{n\rightarrow
\infty}\mathbb{E}\bigg[\Big(L_{\frac{1}{n+p}}(y, I, Y
)-L_{\frac{1}{n}}(y, I, Y)\Big)^{2}\bigg]=0
\end{array}
\end{equation*}
if and only if (\ref{Eq:EY1}) holds. Then $\{L_{\varepsilon}(y, I, Y), \varepsilon>0\}$ is a
Cauchy sequence in $ L^{2}(\Omega, \mathbb{P})$ if and only if
that (\ref{Eq:EY1}) holds. This finishes the proof.
\hfill{$\square$}

\vskip 0.2cm

Now we provided a sufficient and necessary condition for the existence of the
local time of $X$, which complements Theorem 8.1 of Xiao (2009) and
Theorem 3.1 of Jiang and Wang (2009).

\begin{thm}\label{Thm:Existence}
Let $X=\{X(t),\,t\in I\}$ be an $(N,d)$-Gaussian random field defined by (\ref{def:X2}) and
assume that $X_{0}$ has mean zero, continuous covariance function and satisfies
Conditions (C1) and (C2) with index $H \in (0, 1)^N$. Then, for every $x \in \R^d$, 
$L_{\varepsilon}(x, I, X)$ converges in $L^{2}(\Omega, \mathbb{P})$ sense, to a limit 
$L(x, I, X)$ as $\varepsilon\rightarrow 0$ if and only if $\sum_{j=1}^{N}1/H_{j}>d$.
\end{thm}

\pf
By Lemma \ref{Lem:EY2}, we only need to verify that for any $x \in \R^d$,
\begin{equation*}\label{Eq:ETH1}
 \mathcal{M}:=\int_{I^{2}} \exp\bigg(-\frac{\|x\|^{2}\mathbb{E}[(X_{0}(s)-X_{0}(t))^{2}]}{{\rm
detCov}\left(X_0(t), X_0(s)\right)}\bigg)\,\frac{dsdt}{\big[{\rm
detCov}\left(X_0(t), X_0(s)\right)\big]^{\frac d 2}}
\end{equation*}
is finite if and only if $\sum_{j=1}^N\frac1{H_j}>d$.

The sufficiency follows immediately from
Proposition {\ref{prn:prn}}. To prove the necessity, we derive from  (\ref{Eq:prn0}),
(C1) and (C2) that, for any $\varepsilon_0 \in (0, 1)$, there exist constants $c_{11} \ge 1$
 and $c_{12}> 0$  such that $
c_{11}^{-1}\leq {\rm Var}(X_{0}(s))\leq c_{11}$
and
\begin{equation*}\label{Eq:L2con7}
{\rm Var} (X_{0}(t)|X_{0}(s))\ge  c_{12}
\sum_{j=1}^{N}|s_{j}-t_{j}|^{2H_{j}}
\end{equation*}
for all $s, t\in [\varepsilon_{0}, 1]^{N}$. These 
inequalities and (\ref{Eq:CovVar}) imply
\begin{equation}\label{Eq:L2con8}
\arraycolsep=1pt\begin{array}{ll}\displaystyle\frac{\mathbb{E}
[(X_{0}(s)-X_{0}(t))^{2}]}{\det\Gamma_{0,1}(s,t)}\asymp 1
\end{array}
\end{equation}
for all $s, t\in [\varepsilon_{0}, 1]^{N}$.  It follows from
(\ref{Eq:L2con8}) that
\begin{equation*}\label{Eq:L2con9}
\arraycolsep=1pt\begin{array}{ll}
\displaystyle\mathcal{M}&\displaystyle\geq c\, \int_{[\varepsilon_{0},
1 ]^{2N}} \frac{ ds\, dt}{[\det{\rm Cov}(X_{0}(s),
X_{0}(t))]^{d/2} } .
\end{array}
\end{equation*}
From the proof of Proposition \ref{prn:prn} with $\gamma=d$ and
$\lambda=0$ we see that the last integral is infinite if $\sum_{j=1}^N\frac1{H_j}\le d$.
This proves the necessity and hence the theorem.
\hfill{$\square$}

\vskip 0.2cm

In order to study the smoothness of the local times,  
we will make use of the following  lemmas. Lemma \ref{Lem:Hu} 
is from Hu (2001), and Lemma \ref{Lem:ChenYan} is from Chen 
and Yan (2011).
\begin{lem}\label{Lem:Hu}
Let $F\in L^2(\Omega, \mathbb{P})$. Then $F\in D_1$ if and only if
$\Phi_\Theta(1)<\infty.$
\end{lem}

\begin{lem}\label{Lem:ChenYan}
For any $d\in\N$, we have for $x\in[-1,\,1)$,
\begin{equation*}\label{Eq:ChenYan}
\sum_{n=1}^\infty\sum_{
{k_1,\ldots,k_d=0}\atop{k_1+\cdots+k_d=n}}^n\frac{2n(2k_1-1)!!\cdots(2k_d-1)!!}
{(2k_1)!!\cdots(2k_d)!!}x^n\asymp x(1-x)^{-\left(\frac
d2+1\right)}
\end{equation*}
\end{lem}

Recall that the Hermite polynomial of degree $n$ is defined by
\[
H_n(x)=\frac{(-1)^n}{n!}e^{\frac{x^2}{2}}\frac{d^n}{d x^n}
\big(e^{-\frac{x^2}{2}}\big),\quad n\in\Z_+.
\]
It is known that [cf. Nualart (2006)] for any centered Gaussian
random vector $(\xi,\eta)$ with $\E(\xi^2)=\E(\eta^2)=1,$ we have
\begin{equation}\label{Eq:HnHm}
\E\left[H_n(\xi)H_m(\eta)\right]=\left\{
\begin{array}{cc}
0, & m\neq n, \\
\frac1{n!}\left[\E(\xi\eta)\right]^n, & m=n \\
\end{array}
\right.
\end{equation}
and  for all $z\in\C$ and $x\in\R$,
\begin{equation}\label{Eq:Hermite_Exp}
e^{zx-\frac{z^2}{2}}=\sum_{n=0}^\infty z^nH_n(x).
\end{equation}

We will make use of the following lemma.

\begin{lem}\label{Lem:CY2}
Let $Y=\left\{ \big(Y_1(t),\ldots,Y_d(t)\big),\,t\in I \right\}$ 
be an $(N,d)$-Gaussian random field, where $Y_1,\ldots, Y_d$ are 
$d$ independent copies of a centered, real-valued Gaussian random 
field $Y_0$ on $I.$ Suppose that its local time $L(y,I, Y)\in 
L^2(\Omega, \mathbb{P}).$ Then,
\begin{description}
\item {\rm (i)}\,  $L(0, I, Y)\in D_1$ if and only if
\begin{equation}\label{Eq:Smooth}
\int_{I^2}\frac{\big[\E\left(Y_0(s)Y_0(t)\right)\big]^2}{\big[{\rm
detCov}\left(Y_0(t), Y_0(s)\right)\big]^{\frac d 2+1}}\,ds\,dt<\infty.
\end{equation}
\item {\rm (ii)}\,  If  (\ref{Eq:Smooth}) holds, then $L(y, I, Y)\in D_1$ for every
$y \in \R^d \backslash\{0\}$.
\end{description}
\end{lem}
\pf The proof is similar to that of Lemma 3.2 in Chen and Yan (2011),
see also Hu (2001). Let $L_\eps(y,I,Y)$ be as in (\ref{Def:selfv}) (by replacing $X$ by $Y$).
Thanks to (\ref{Def:selfv}) and (\ref{Eq:Hermite_Exp}), we can write
\begin{eqnarray}\label{Eq:L_eps1}
L_\eps(y,I, Y)
&=&\frac1{(2\pi)^d}\int_I\int_{\R^d}e^{-i\l\xi,y\r} \exp\bigg( i\l\xi,Y(s)\r-\eps\frac{\|\xi\|^2}{2}\bigg)\,d\xi ds \nonumber\\
&=&\frac1{(2\pi)^d}\int_I\int_{\R^d}e^{-i\l\xi,y\r} \exp\Big(-\frac12\left(\E(Y_0^2(s))+\eps\right)\|\xi\|^2\Big)\nonumber \\
&\ &\quad\times \sum_{n=0}^\infty i^n\left(\E(Y_0^2(s))\|\xi\|^2\right)^{n/2}H_n\bigg(\frac{\l\xi, Y(s)\r}
{\sqrt{\E(Y_0^2(s))\|\xi\|^2}}\bigg)\,d\xi ds\nonumber\\
&:=&\sum_{n=0}^\infty F_n^{y,\eps}.
\end{eqnarray}
Denote $\Phi_{\Theta_{y, \eps}}(u)=\E\left(|\Gamma_{\sqrt u}  L_\eps(y,I,Y)|^2\right)$ and
$\Phi_{\Theta_y}(u)=\E\left(|\Gamma_{\sqrt u} L(y, I, Y)|^2\right).$  
Also, for simplicity of notation, let $a^2 = \E(Y_0^2(s))+\eps$ and $b^2 = \E(Y_0^2(t))+\eps$.  
It follows from (\ref{Eq:L_eps1}) and (\ref{Eq:HnHm}) that
\begin{equation}\label{Eq:Phi1}
 \begin{split}
&\Phi_{\Theta_{y,\eps}}(1)=\sum_{n=0}^\infty n\E\left(|F_n^{y,\eps}|^2\right)\\
&=\sum_{n=0}^\infty\frac{n}{(2\pi)^{2d}}\E\Bigg[\int_{I^2}\int_{\R^{2d}}
e^{-i\l\xi -\eta ,y\r }\left[\E(Y_0^2(s))\E(Y_0^2(t))\|\xi\|^2\|\eta\|^2\right]^{n/2}\\
&\qquad\times\exp\Big(-\frac12\big[a^2\|\xi\|^2+ b^2\|\eta\|^2\big]\Big)\\
&\qquad\times H_n\bigg(\frac{\l\xi, Y(s)\r}{\sqrt{\E(Y_0^2(s))\|\xi\|^2}}\bigg) 
H_n\bigg(\frac{\l\eta, Y(t)\r}{\sqrt{\E(Y_0^2(t))\|\eta\|^2}}\bigg) \,d\xi d\eta ds dt \Bigg]\\
&= \sum_{n=1}^\infty\frac{1}{(2\pi)^{2d}(n-1)!}\int_{I^2} \left[\E\left(Y_0(s)Y_0(t)\right)\right]^n\, dsdt\\
 &\qquad \times \int_{\R^{2d}} e^{-i\l\xi -\eta ,y\r} \, \l\xi,\eta\r^n
\exp\Big(-\frac12\big[a^2 \|\xi\|^2+ b^2|\eta\|^2\big]\Big)\,d\xi d\eta.
\end{split}
\end{equation}

If $y = 0$, then the integrals in (\ref{Eq:Phi1}) become 0 for all odd numbers $n$. Hence
\begin{equation}\label{Eq:Phi1b}
\begin{split}
\Phi_{\Theta_{0, \eps}}(1)&= \sum_{n=1}^\infty \frac{1}{(2\pi)^{2d}(2n-1)!}
\int_{I^2}\left[\E\left(Y_0(s)Y_0(t)\right)\right]^{2n}\,ds dt  \\
&\quad\times\int_{\R^{2d}} \l\xi,\eta\r^{2n} \exp\Big(-\frac12\big[a^2\|\xi\|^2+
b^2\|\eta\|^2\big]\Big)\,d\xi d\eta.
\end{split}
\end{equation}
By using the fact that for $k\in\Z_+$, $\gamma>0,$
\[
\int_\R v^{2k}\exp\Big(-\frac{\gamma v^2}2\Big)dv=\sqrt{2\pi}(2k-1)!!\gamma^{-(k+1/2)}
\]
and the same argument as in Chen and Yan (2011, p. 1010), we obtain
\begin{equation}\label{Eq:Phi1_1}
\begin{split}
& \frac{1} {(2\pi)^d }\int_{\R^{2d}}\frac{ \l\xi,\eta\r^{2n} }{(2n-1)!}
\exp\Big(-\frac12\big[a^2 \|\xi\|^2+ b^2 \|\eta\|^2\big]\Big)\,d\xi d\eta\\
&= \sum_{{k_1,\ldots,k_d=0}\atop{k_1+\cdots+k_d=n}}^n\frac{2n(2k_1-1)!!\cdots(2k_d-1)!!}
{(2k_1)!!\cdots(2k_d)!!\big[\left(\E(Y_0^2(s))+\eps\right)
\left(\E(Y_0^2(t))+\eps\right)\big]^{n+\frac d2 }}.
\end{split}
\end{equation}
[This can be verified by using induction.] Combining (\ref{Eq:Phi1b}) and (\ref{Eq:Phi1_1}), 
and applying Lemma \ref{Lem:ChenYan} and the monotone convergence theorem, we obtain
\begin{equation}\label{Eq:Phi1_2}
\begin{split}
& \sum_{n=1}^\infty\frac{1}{(2n-1)!}\int_{I^2}\left[\E\left(Y_0(s)Y_0(t)\right)\right]^{2n}\,ds dt
\int_{\R^{2d}} \l\xi,\eta\r^{2n} \\
& \ \qquad \quad \times \exp\Big(-\frac12\big[a^2\|\xi\|^2+
b^2\|\eta\|^2\big]\Big)\,d\xi
d\eta \\
&  \asymp
\int_{I^2}\frac{\left[\E\left(Y_0(s)Y_0(t)\right)\right]^{2}}
{\big\{a^2 b^2 -\left[\E\left(Y_0(s)Y_0(t)\right)\right]^{2}\big\}^{\frac d 2+1}}\,dsdt  \\
& \stackrel{\mbox{ as } \eps\to 0}{\longrightarrow}  \int_{I^2}\frac{\left[\E\left(Y_0(s)Y_0(t)\right)\right]^{2}}
{\big\{\E(Y_0^2(s))\E(Y_0^2(t))-\left[\E\left(Y_0(s)Y_0(t)\right)\right]^{2}\big\}^{\frac d 2+1}}\,dsdt   \\
& =\int_{I^2}\frac{\big[\E\left(Y_0(t)Y_0(s)\right)\big]^2}{\big[{\rm
detcov}\left(Y_0(s), Y_0(t)\right)\big]^{\frac d2+1}}\,ds\,dt,
\end{split}
\end{equation}
which proves Part (i) of Lemma \ref{Lem:CY2}, thanks to Lemma \ref{Lem:Hu}.

Now we prove Part (ii) of the lemma. Notice that for $y \in \R^d \backslash\{0\}$, 
it does not seem easy to compute the integrals in the last equality of (\ref{Eq:Phi1}) 
explicitly. So we turn to the following upper bound.
\begin{equation}\label{Eq:Phi1_3}
\begin{split}
\Phi_{\Theta_{y,\eps}}(1)& \le \sum_{n=1}^\infty\frac{1}{(2\pi)^{2d}(n-1)!}
\int_{I^2} \big|\E\left(Y_0(s)Y_0(t)\right)\big|^n\, dsdt\\
&\qquad  \times  \int_{\R^{2d}}   \big| \l\xi,\eta\r \big|^n
\exp\Big(-\frac12\big[a^2\|\xi\|^2+ b^2
 \|\eta\|^2\big]\Big)\,d\xi d\eta.
\end{split}
\end{equation}
The sum over even integers in (\ref{Eq:Phi1_3}) is the same as in (\ref{Eq:Phi1_2}). So we only need to
consider the terms over odd integers. For this purpose, let
\begin{equation}\label{Eq:Jk}
\begin{split}
J_{2n+1} &= \frac 1 {(2n)!} \int_{I^2} \left| \E\left(Y_0(s)Y_0(t)\right)\right|^{2n+1}  \, dsdt  \\
&\qquad \times  \int_{\R^{2d}}
 |\l\xi,\eta\r |^{2n+1} \exp\Big(-\frac12\big[a^2\|\xi\|^2+
b^2\|\eta\|^2\big]\Big) \,d\xi d\eta.
\end{split}
\end{equation}
By using the Cauchy-Schwarz inequality and the elementary inequality
$$
x e^{-\frac \beta {2n}   x^2}  \le \sqrt{\frac n {e \beta}}\, , \qquad  \forall \beta > 0\ \hbox{ and } \,
x > 0,
$$
we have
\[
|\l\xi,\eta\r | \, e^{-\frac1 {2n}\big[a^2\|\xi\|^2+
b^2\|\eta\|^2\big]} \le \frac {n}
{e\, ab}.
\]
Plugging this into (\ref{Eq:Jk}) yields
\begin{equation}\label{Eq:Jk2}
\begin{split}
J_{2n+1} &\le   \frac {1} {2e (2n-1)!}  \int_{I^2} \left| \E\left(Y_0(s)Y_0(t)\right)\right|^{2n} \, dsdt \\
&\qquad \times  \int_{\R^{2d}}
 |\l\xi,\eta\r |^{2n} \, \exp\Big(-\frac{n-1} {2n} \big[a^2\|\xi\|^2+
b^2\|\eta\|^2\big]\Big) \,d\xi d\eta.
\end{split}
\end{equation}
The same argument for \eqref{Eq:Phi1_1} gives
\begin{equation}\label{Eq:Jk3}
\begin{split}
& \int_{\R^{2d}}  \frac { |\l\xi,\eta\r |^{2n}} {(2n-1)!}
 \exp\bigg(-\frac{n-1}{2n} \big[a^2\|\xi\|^2+
b^2\|\eta\|^2\big]\bigg) \,d\xi d\eta  \\
&=\frac{(2\pi)^d}{ (1 - n^{-1})^{n + \frac d 2}}\sum_{{k_1,\ldots,k_d=0}
\atop{k_1+\cdots+k_d=n}}^n\frac{2n(2k_1-1)!!\cdots(2k_d-1)!!}
{(2k_1)!!\cdots(2k_d)!!\big[a^2
b^2\big]^{n+\frac d2}}.
\end{split}
\end{equation}
Combining  (\ref{Eq:Jk})--(\ref{Eq:Jk3})  with (\ref{Eq:Phi1_3}), and using 
the same argument as in (\ref{Eq:Phi1_2}), we derive $\Phi_{\Theta_{y}}(1) < \infty$ 
under (\ref{Eq:Smooth}). 
This finishes the proof of Part (ii).
\hfill{$\square$}

\vskip 0.2cm

The following is the main theorem of this section.

\begin{thm}\label{Thm:Smooth}
Let $X=\{X(t),\,t\in I\}$ be an $(N,d)$-Gaussian  field
defined by (\ref{def:X2}) and assume that $X_{0}$ satisfies
(C1) and (C2) with index $H \in (0, 1)^N$. Then the following statements hold:
\begin{description}
\item {\rm (i)}\,  $L(0,I, X)\in D_1$ if and only if $\sum_{j=1}^N\frac1{H_j}>d+2.$
\item {\rm (ii)}\,  If $\sum_{j=1}^N\frac1{H_j}>d+2$, then $L(x,I, X)\in D_1$  for every
$x \in \R^d \backslash \{0\}$.
\end{description}
\end{thm}

\pf By Theorem \ref{Thm:Existence}  and Lemma \ref{Lem:CY2}, it is sufficient for us to
verify that
\begin{equation}\label{Eq:Smooth1}
\int_{I^2}\frac{\big[\E\left(X_0(s)X_0(t)\right)\big]^2}{\big[{\rm
detCov}\left(X_0(t), X_0(s)\right)\big]^{\frac
d2+1}}\,ds\,dt<\infty
\end{equation}
if and only if  $\sum_{j=1}^N\frac1{H_j}>d+2$. This  follows from Proposition \ref{prn:prn} with
$\gamma=d+2$ and $\lambda=2$ immediately. \hfill{$\square$}

\vskip 0.2cm
\begin{rem}\,
As we mentioned in Introduction,  the class of Gaussian random fields that satisfy (C1) and (C2) is large,
including fractional Brownian sheets and the solutions to stochastic heat equation driven by various 
space-time  Gaussian noises. In particular,  Theorem \ref{Thm:Smooth}  recovers Theorem 11 in Eddahbi, 
et al. (2005)  and Theorem 2.1 with $\alpha = 1$ in Eddahbi and Vives (2003).
\end{rem}

In the following we apply Theorems \ref{Thm:Existence} and
\ref{Thm:Smooth} to study the collision and intersection local 
times of independent Gaussian fields. Theorems  
\ref{Thm:SmoothCol} and \ref{Thm:SmoothInt} below generalize 
the results of Yan et al. (2009), Yan and Shen (2010) and Chen and
Yan (2011) for fractional Brownian motion and related Gaussian processes.

\vskip 0.4cm

\noindent{\bbb 2.2\quad Smoothness of the collision local time}\\


Given $H=(H_{1}, \ldots, H_{N})\in (0,1)^{N}$ and $K=(K_{1}, \ldots,
K_{N})\in (0,1)^{N}$, let
$X^{H}=\{X^{H}(s), s\in \mathbb{R}^{N}\}$ and $X^{K}=\{X^{K}(t),
t\in \mathbb{R}^{N}\}$ be two independent Gaussian random fields
with values in $\mathbb{R}^{d}$ as defined in (\ref{def:X2}). We
assume that the associate real-valued random fields $X^{H}_0$ and
$X^{K}_0$ satisfy Conditions (C1) and (C2) on interval $I
\subseteq \R^{N}$ respectively with indices $H$ and  with indices $K$.

The collision local time of $X^H$ and $X^K$ on $I$ is
formally defined by
\begin{equation}\label{Def:Collision}
L_C(X^H,X^K):=\int_I\delta\left(X^H(s)-X^K(s)\right)\,ds.
\end{equation}

\begin{thm}\label{Thm:SmoothCol} Let $L_C\left(X^H,X^K\right)$ be the
collision local time of $X^H$ and $X^K$ as above. Then
\begin{itemize}
\item[(i)]\, $L_C\left(X^H,X^K\right)\in L^2(\Omega,\mathbb{P})$
if and only if $\sum_{j=1}^N\frac1{H_j \wedge K_j}>d.$
\item[(ii)]\, $L_C\left(X^H,X^K\right)\in D_1$ if and only if
$\sum_{j=1}^N\frac1{H_j\wedge K_j}>d+2.$
\end{itemize}
\end{thm}

\pf Consider the $(N, d)$ Gaussian field $Z =
\{Z(t), t \in I\}$ defined by 
$$Z(t)\equiv X^H(t)-X^K(t), \qquad \forall \ t\in I.$$
Then the collision local time of $X^H$ and $X^K$ on $I$ is nothing but
$L(0, I, Z)$, the  local time of $Z$ on $I$ at $x=0$. Hence Theorem
\ref{Thm:SmoothCol} follows from Theorems \ref{Thm:Existence}
and \ref{Thm:Smooth} once we verify that the real valued
Gaussian field $Z_0(t) = X_0 (t)- Y_0(t)$ satisfies
(C1) and (C2) in the interval $I$ with
indices $(H_1 \wedge K_1, \ldots, H_N \wedge K_N)\in(0,1)^{N}$.
%

Since it is easy to show that $Z_{0} $ satisfies (C1), we verify (C2) only.
By the definition of conditional variance and independence of $X^H$ and
$X^K$, we have
\begin{equation*}\label{Eq:Z3}
\begin{split}
&{\rm Var}\left(Z_0(t)|Z_0(s)\right)
=\inf_{a\in\R}\left\{\E\left[\big(X^H_0(t)-aX^H_0(s)\big)^2+ \big(X^K_0(t)-aX^K_0(s)\big)^2\right]\right\}\\
& \qquad \quad \ge \inf_{a\in\R}\E\left[\big(X^H_0(t)-aX^H_0(s)\big)^2\right]+
\inf_{b\in\R}\E\left[\big(X^K_0(t)-bX^K_0(s)\big)^2\right]\\
& \qquad \quad={\rm Var}\left(X_0^H(t)|X_0^H(s)\right)+{\rm Var}\left(X_0^K(t)|X_0^K(s)\right)\\
& \qquad \quad\ge c\, \bigg(\sum_{j=1}^N\min\big\{|s_j-t_j|^{2H_j},
t_j^{2H_j}\big\}+ \sum_{j=1}^N\min\big\{|s_j-t_j|^{2K_j},
t_j^{2K_j}\big\}\bigg)
\\
& \qquad \quad\ge c \sum_{j=1}^N\min\big\{|s_j-t_j|^{2(H_j\wedge
K_j)}, t_j^{2(H_j\wedge K_j)}\big\},\quad \forall s,\,t\in I,
\end{split}
\end{equation*}
for some constant $c >0.$  This verifies  that
$Z_{0}$ satisfies  Condition (C2).  \hfill{$\square$}

\vskip 0.4cm

\noindent{\bbb 2.3\quad Smoothness of the intersection local time }\\[0.1cm]


Let $H=(H_{1}, \ldots, H_{N_{1}})\in (0,1)^{N_{1}}$ and $K=(K_{1},
\ldots, K_{N_{2}})\in (0,1)^{N_{2}}$ be two constant vectors. Let
$X^{H}=\{X^{H}(s), s\in \mathbb{R}^{N_{1}}\}$ and
$X^{K}=\{X^{K}(t), t\in \mathbb{R}^{N_{2}}\}$ be two independent
Gaussian random fields with values in $\mathbb{R}^{d}$ as defined
in (\ref{def:X2}). We assume that the associate real-valued random
fields $X^{H}_0$ and $X^{K}_0$ satisfy Conditions (C1) and (C2)
respectively on interval $I_1 \subseteq \R^{N_1}$ with indices
$H=(H_{1}, \ldots, H_{N_{1}})$ and on $I_2 \subseteq \R^{N_2}$
with indices $K=(K_{1}, \ldots, K_{N_{2}})$.
Then the intersection local time of $X^H$ and $X^K$ is formally defined by
\begin{equation}\label{Def:Intersection}
L_I(X^H,X^K):=\int_{I_{N_1}\times
I_{N_2}}\delta\left(X^H(s)-X^K(t)\right)\,ds\,dt.
\end{equation}


\begin{thm}\label{Thm:SmoothInt}
Let $L_I\left(X^H,X^K\right)$ be the intersection
local time of $X^H$ and $X^K$ as above. Then
\begin{description}
\item {\rm (i)}\, $L_I\left(X^H,X^K\right)\in L^2 (\Omega,\mathbb{P})$ 
if and only if $\sum_{j=1}^{N_1}\frac1{H_j}
+\sum_{j=1}^{N_2}\frac{1}{K_j}>d.$
\item {\rm (ii)}\, $L_I\left(X^H,X^K\right)\in D_1$
if and only if $\sum_{j=1}^{N_1}\frac1{H_j}
+\sum_{j=1}^{N_2}\frac{1}{K_j}>d+2.$
\end{description}
\end{thm}

\pf  Let $U=\{U(s,t),\,(s,t)\in I_{N_1}\times I_{N_2}\}$ be the
$(N_1+N_2,d)$-Gaussian random field with mean $0$ defined by
$$ U(s,t)=X^H(s)-X^K(t), \quad \forall \, s\in I_{N_1}, \, t\in I_{N_2}. $$
Clearly, the intersection local time of $X^H$ and $X^K$ is nothing
but $L(0, I_{N_1} \times I_{N_2}, U),$ the local time of $U$ on
$I_{N_1} \times I_{N_2}$ at $x = 0$.
One can verify  that  the Gaussian random field
 $U_{0}(s,t)=X^H_0(s)-X^K_0(t) $   satisfies Conditions (C1)
and (C2)  on the interval $I_{N_1}\times I_{N_2}$ with indices
$(H_1,\ldots,H_{N_1},K_1,\ldots,K_{N_2})\in (0,1)^{N_{1}+N_2}$.
Therefore, the conclusions follow from Theorems \ref{Thm:Existence} and
\ref{Thm:Smooth}. \hfill{$\square$}

\noindent\\[2mm]

\setcounter{equation}{0}
\setcounter{section}{3}
\setcounter{thm}{0}

\noindent{\bbb 3\quad Self-intersection local times}\\[0.1cm]

In this section, we study the existence and smoothness of self-intersection local times
of an $(N, d)$-Gaussian random field $X =\{X(t), t\in \mathbb{R}^{N}\}$ as in (\ref{def:X2}).
These problems are more involved than the collision or intersection local times of
independent Gaussian random fields, due to complexity of dependence structures
of $X$. For earlier results for the Brownian sheet,
fractional Brownian motion and related self-similar Gaussian processes, we refer to
Imkeller and Weise (1995, 1999), Hu (2001), Jiang and Wang (2009). Their methods rely on
special properties of the Brownian sheet or fractional Brownian motion. Our approach
below is based on a weak form of local nondeterminism and is more general.

For any two compact intervals $I,  J\subseteq \R^{N}$, the self-intersection
local times of $X =\{X(t), t\in \mathbb{R}^{N}\}$ on $I$ and $J$ is formally defined by
\begin{equation}\label{Def:self}
\arraycolsep=1pt\begin{array}{ll} L_{S}(X, I\times
J)&\displaystyle=\int_{I\times J}\delta(X(s)-X(t))dsdt.
\end{array}
\end{equation}

Define a $(2N, d)$-Gaussian random field $V=\{V(s,t), (s,t)\in \mathbb{R}^{2N}\}$
by
\begin{equation}\label{def:V(s t)}
V(s,t):= X(s)-X(t), \ s, t\in \mathbb{R}^{N}.
\end{equation}
Then the self-intersection local time of $X$ is $L(0, I \times J, V), $
the local time of $V$  on $I \times J$ at $x = 0$.

Under the condition that $X_0$ satisfies Conditions (C1) and (C2) on both
intervals $I$ and $J$, the Gaussian field $V_{0}(s,t) = X_0(s) - X_0(t)$
may not satisfy the corresponding (C2) on $I\times J$. Therefore,
we can not apply Theorems \ref{Thm:Existence} and \ref{Thm:Smooth} directly.
To overcome this difficulty, we will make use of
the following condition:
\begin{itemize}
\item[(C3)]\ There exists a positive
constant $c_{12}$ such that for all $u, t^1, t^2, t^3\in [0, 1
]^{N},$
\begin{equation}\label{Eq:LND}
{\rm Var}\left(X_0(u)|X_0(t^1), X_0(t^2),
X_0(t^3)\right)\ge c_{12}\sum_{j=1}^N\min_{0\le k\le
3}|u_j-t_j^k|^{2H_j},
\end{equation}
\ \ where $t_{j}^{0}=0$, \ $j=1, 2, \ldots, N$.
\end{itemize}

Clearly, Condition (C2) is a special case of Condition (C3). It is known that
multiparameter fractional Brownian motion and fractional Brownian sheets
satisfy Conditions (C1) and (C3); see Pitt (1978) and Wu and Xiao (2007).
More examples can be found in Xiao (2009).

For two compact intervals $I,\,J \subseteq[0,\,1]^N$, we call them {\it separated} if
\begin{equation}\label{Eq:Sprt}
\inf_{s\in I,\,t\in J}|s_j-t_j|>0\,\,\mbox{ for some }\,\,j=1,\,2,\ldots,\, N.
\end{equation}
Let $S\subseteq\{1,\ldots,N\}$ be the collection of all $j$'s that satisfy (\ref{Eq:Sprt})
and let  $S^c=\{1,\ldots,N\}\setminus S.$ Because $I$ and $J$ are compact, there exists
$\varepsilon_0>0$ such that
\begin{equation}\label{Eq:Sprt2}
\inf_{s\in I,\,t\in J}|s_j-t_j| \ge\varepsilon_0\,\,\mbox{ for }\,\,j\in S.
\end{equation}

We further call $I$ and $J$  {\it partially separated} if both $S$ and $S^c$ are nonempty,
{\it well separated} if $S^c=\emptyset$, and {\it not separated} if $S=\emptyset.$
Clearly, $I$ and $J$ are not separated iff $I\cap J\neq\emptyset.$

Similarly to Imkeller and Weisz (1999) for the Brownian sheet, we consider the
self-intersection local times of $X$ on $I$ and $J$ by distinguishing
three cases:
\begin{description}
\item {\rm Case (i)}\ \  $I,\,J\subseteq[0,\,1]^N$ are well separated.
\item {\rm Case (ii)}\   $I,\,J\subseteq[0,\,1]^N$ are partially separated.
\item {\rm Case (iii)}   $I,\,J\subseteq[0,\,1]^N$ are not separated.
\end{description}

In Case (i), we have the following theorem.
\begin{thm}\label{Thm:SmoothSelfO}
Let $X=\{X(t), t\in \mathbb{R}^{N}\}$ be an $(N, d)$-Gaussian random field
defined by (\ref{def:X2}) with $X_0$ satisfying Conditions (C1) and (C3)
and let $L_S\left(X,I\times J\right)$ be the self-intersection local time
of $X$ on $I$ and $J$. If $I$ and $J$ are well separated, then the following
statements hold:
\begin{description}
\item {\rm (i)}\ \ $L_S\left(X,I\times J\right)\in L^2(\Omega, \mathbb{P})$
if and only if $2\sum_{j=1}^N\frac1{H_j}>d.$
\item {\rm (ii)}\ $L_S\left(X,I\times J\right)\in D_1$
if and only if $2\sum_{j=1}^N\frac1{H_j}>d+2.$
\end{description}
\end{thm}

\pf
Since the Gaussian field $X_0$ satisfies (C1) on $I$ and $J$, we see
that for any $ (s,t),(s', t')\in I\times J$,
\begin{equation}\label{Eq:V2}
\begin{split}
\E\left[\big(V_0(s,t)-V_0(s',t')\big)^2\right]&\le
c\, \bigg[\sum_{j=1}^{N}|s_j-s'_j|^{2H_j}+\sum_{j=1}^{N}|t_j-t'_j|^{2H_j}\bigg].
\end{split}
\end{equation}
Thus the Gaussian field $V_{0}(s,t)=X_0^H(s)-X_0^H(t)$ satisfies (C1)
 on $I\times J$ with indices
 $(H_1,\ldots,H_N,\,H_1,\ldots,H_N)$ $\in(0,\,1)^{2N}.$ To verify that
 $V_0$ also satisfies (C2), we see that (C3) implies
\begin{equation*}\label{Eq:V3}
\begin{split}
{\rm Var}\left(V_0(s,t)|V_0(s',t')\right)
&\ge {\rm Var}\left(X_0(t)|X_0 (s),\,X_0(s'),\,X_0(t')\right)\\
&\ge c_{12} \sum_{j=1}^{N}\min\big\{|t_j-s_j|^{2H_j},
|t_j-s'_j|^{2H_j}, |t_j-t'_j|^{2H_j},  t_j^{2H_j}\big\}\\
&\ge c_{13}\, \sum_{j=1}^{N}\min\big\{
|t_j-t'_j|^{2H_j},  t_j^{2H_j}\big\}
\end{split}
\end{equation*}
thanks to the fact that $|t_j-s_j|\ge \varepsilon_{0}$ and
$|t_j-s'_j|\ge \varepsilon_{0}$.  Here the constant $c_{13}$ depends on $\varepsilon_0$.
By the same token, we have
\begin{equation*}\label{Eq:V4}
{\rm Var}\left(V_0(s,t)|V_0(s',t')\right)\ge
c_{13}\sum_{j=1}^{N}\min\big\{|s_j-s'_j|^{2H_j},
s_j^{2H_j}\big\}.
\end{equation*}
Adding up these two inequalities shows
\begin{equation*}\label{Eq:V5}
\begin{split}
&{\rm Var}\left(V_0(s,t)|V_0(s',t')\right)\\
&\ge
\frac{c_{13}} 2\, \bigg[\sum_{j=1}^{N}\min\big\{|s_j-s'_j|^{2H_j},
s_j^{2H_j}\big\}+\sum_{j=1}^{N}\min\big\{ |t_j-t'_j|^{2H_j},
t_j^{2H_j}\big\}\bigg].
\end{split}
\end{equation*}
This proves that $V_0$ satisfies (C2)
on $I\times J$ with $(H_1,\ldots,H_N,\,H_1,\ldots,H_N)$ $\in(0,\,1)^{2N}.$
Therefore, the conclusions follow from  Theorems \ref{Thm:Existence} and
\ref{Thm:Smooth}. \hfill{$\square$}

\vskip 0.2cm

Now we  consider Case (ii), e.g.  the two compact intervals
$I$ and $J$  are partially separated. In this case, both $S$ and $S^c$ are
nonempty sets. For concreteness,  we may assume that
$I=[a,\,a+\l h\r],\,\,J=[b,\,b+\l h\r]$, where $b_j > a_j +h$ for $ j \in S$ and $a_j = b_j $
 for $ j \in S^c$. Then (\ref{Eq:Sprt2}) holds with $\varepsilon_0
 = \min \{b_j - a_j - h, {j \in S}\}$. Note that, when $X$ is the $(N, d)$
Brownian sheet, the existence condition in (i) in the following theorem
coincides with that
in Theorem 3 of Imkeller and Weisz (1995, 1999).

\begin{thm}\label{Thm:SmoothSelfO2}
Let $X=\{X(t), t\in \mathbb{R}^{N}\}$ be an $(N, d)$-Gaussian random field
as in Theorem \ref{Thm:SmoothSelfO}. Let $I$ and $J$ be partially separated
as described above.
Then the following statements hold:
\begin{description}
\item {\rm (i)}\  $L_S\left(X,I\times J\right)\in L^2(\Omega, \mathbb{P})$ if  \,
$2\sum_{j\in S}\frac1{H_j}+\sum_{j\in S^c}\frac1{H_j}>d.$
\item {\rm (ii)}\  $L_S\left(X,I\times J\right)\notin L^2(\Omega, \mathbb{P})$
if\, $2\sum_{j=1}^N\frac1{H_j}\le d.$
\item {\rm (iii)}\  $L_S\left(X,I\times J\right)\in D_1$
if\, $2\sum_{j\in S}\frac1{H_j}+\sum_{j\in S^c}\frac1{H_j}>d+2.$
\item {\rm (iv)}\ $L_S\left(X,I\times J\right)\notin D_1$
if\,  $2\sum_{j=1}^{N}\frac1{H_j}\leq d+2.$
\end{description}
\end{thm}

\pf We prove Part (i) at first. By Lemma \ref{Lem:EY2},  we only need to
prove that if $2\sum_{j\in S}\frac1{H_j}+\sum_{j\in S^c} \frac1{H_j}$ $>d$ then
\begin{equation}\label{Eq:2S1}
\arraycolsep=1pt\begin{array}{ll} \mathcal{J}:=\displaystyle\int_{
(I\times J)^{2}}\frac{dsdtds'dt'}{[\det{\rm Cov}(V_{0}(s, t), V_{0}(s',
t'))]^{d/2} }<\infty.
\end{array}
\end{equation}
By the definition of conditional variance
and (C3), we see that for any $ (s,t),\,(s',t')\in I\times J$,
\begin{equation}\label{Eq:2S2}
\begin{split}
&{\rm Var}\left(V_0(s,t)|V_0(s',t')\right)
\ge {\rm Var}\left(X_0(t)|X_0(s),\,X_0(s'),\,X_0(t')\right)\\
&  \qquad\ge c_{12}\sum_{j=1}^{N}\min\big\{|t_j-s_j|^{2H_j},
|t_j-s'_j|^{2H_j}, |t_j-t'_j|^{2H_j}, t_j^{2H_j}\big\}\\
& \qquad \ge c_{14}\bigg(\sum_{j\in S}\min\big\{|t_j-t'_j|^{2H_j}, t_j^{2H_j}\big\}\\
&\qquad \quad + \sum_{j\in S^c}\min\big\{|t_j-s_j|^{2H_j},
|t_j-s'_j|^{2H_j}, |t_j-t'_j|^{2H_j}, t_j^{2H_j}\big\}\bigg),
\end{split}
\end{equation}
thanks to the fact that if $j\in S$, then $|t_j-s_j|\ge \varepsilon_0$
and $|t_j-s'_j|\ge \varepsilon_0.$ By the same token, we have
\begin{equation}\label{Eq:2S3}
\begin{split}
&{\rm Var}\left(V_0(s,t) \big|V_0(s',t')\right)
\ge c_{14} \bigg(\sum_{j\in S}\min\big\{|s_j-s'_j|^{2H_j}, s_j^{2H_j}\big\}\\
&\qquad \qquad +\sum_{j\in S^c}\min\big\{|s_j-t_j|^{2H_j},
|s_j-s'_j|^{2H_j}, |s_j-t'_j|^{2H_j}, \, s_j^{2H_j}\big\}\bigg).
\end{split}
\end{equation}
Combining (\ref{Eq:2S2}) and (\ref{Eq:2S3}), we  have
\begin{equation}\label{Eq:2S4}
\begin{split}
&{\rm Var}\left(V_0(s,t) \big|V_0(s',t')\right)\\
&\ge c_{14} \bigg[\sum_{j\in S} \Big(\min\big\{|t_j-t'_j|^{2H_j}, t_j^{2H_j}\big\} +
\min\big\{|s_j-s'_j|^{2H_j}, s_j^{2H_j}\big\}\Big) \\
&\qquad+\sum_{j\in S^c}\min\big\{|s_j-t_j|^{2H_j},
|s_j-s'_j|^{2H_j}, |s_j-t'_j|^{2H_j}, s_j^{2H_j}\big\}\bigg].
\end{split}
\end{equation}
Note that, in (\ref{Eq:2S4}), only one sum over $S^c$ in (\ref{Eq:2S2})
and (\ref{Eq:2S3}) is kept. This is due to the fact that, when integrating
$ds_j$ for $j \in S^c$, all the other variables,  $s_j',\,t_j$ and $ t_j'$,
will disappear [see Lemma 2.2]. This situation is different
from the case when we integrate $ds_j$ for $j \in S$.

Since $I$ and $J$ are partially separated (i.e., $S \ne \emptyset$), we have
\begin{equation}\label{Eq:2S5}
{\rm Var}\left(V_0(s',t')\right)=\E\left[\big(X_0(s')-X_0(t')\big)^2\right]
 \asymp 1,\quad \forall \ s'\in I,\, t'\in J.
\end{equation}
It follows from (\ref{Eq:2S4}) and (\ref{Eq:2S5}) that the integral $\mathcal{J}$ in
(\ref{Eq:2S1}) is at most
\begin{equation}\label{Eq:2S6}
\arraycolsep=1pt\begin{array}{ll}  &\displaystyle
\int_{(I\times J)^{2}}\bigg[\sum_{j\in S}\Big(\min\big\{|t_j-t'_j|^{H_j}, t_j^{H_j}\big\}+
\min\big\{|s_j-s'_j|^{H_j}, s_j^{H_j}\big\}\Big)\\
&\displaystyle \quad + \sum_{j\in S^c}\min\big\{|s_j-t_j|^{H_j},
|s_j-s'_j|^{H_j}, |s_j-t'_j|^{H_j}, s_j^{H_j}\big\}\bigg]^{-d}\, dsdtds'dt'.
\end{array}
\end{equation}
Similarly to the argument in the proofs of (\ref{Eq:prn2}) and (\ref{Eq:suf1}),
we integrate iteratively and apply Lemmas 2.1-2.3 to show that the integral in
(\ref{Eq:2S6}) is finite if $2\sum_{j\in S=1}\frac1{H_j} +\sum_{j\in S^c}\frac1{H_j}>d$.
This proves the sufficiency in Part (i).

Next we prove Part (ii). For any $(s,t), (s',t')\in  I\times J$,
Condition (C1) implies that
\begin{equation}\label{Eq:2S7}
\begin{split}
{\rm Var}\left(V_0(s,t)|V_0(s',t')\right)
&\le\E\left[\big(X_0(s)-X_0(t)-X_0(s')+X_0(t')\big)^2\right]\\
& \le c\, \sum_{j=1}^{N}\left(|s_j-s'_j|^{2H_j}+|t_j-t'_j|^{2H_j}\right).
\end{split}
\end{equation}
It follows from (\ref{Eq:2S5}) and (\ref{Eq:2S7}) that
\begin{equation}\label{Eq:2S8}
\det{\rm Cov}(V_{0}(s, t), V_{0}(s', t'))\le
c\, \sum_{j=1}^{N}\left(|s_j-s'_j|^{2H_j}+|t_j-t'_j|^{2H_j}\right).
\end{equation}
This implies, by using Lemma 2.1 repeatedly, that the integral $\mathcal{J}$
in (\ref{Eq:2S1}) is infinite if $2\sum_{j=1}^N\frac1{H_j} \le d$.

In order to prove Part (iii), by Lemma \ref{Lem:CY2}, it suffices to show that,
if  $2\sum_{j\in S}\frac1{H_j}+\sum_{j\in S^c}\frac1{H_j}>d+2$, then
\begin{equation}\label{Eq:2P1}
\mathcal{K}=\displaystyle\int_{(I\times J)^2}\frac{\big[\E\left(V_{0}(s, t)
V_{0}(s', t')\right)\big]^2}{[\det{\rm Cov}(V_{0}(s, t), V_{0}(s',
t'))]^{\frac{d+2}{2}} }\,ds\,dt\,ds'\, dt'<\infty.
\end{equation}

For any $(s,t), (s',t')\in I\times J$, we use the Cauchy-Schwarz
inequality and  (C1) again to show that
\begin{equation}\label{Eq:2P2}
\begin{split}
& \big[\mathbb{E}(V_{0}(s, t) V_{0}(s', t'))\big]^{2}\le
\mathbb{E}\big[V^{2}_{0}(s, t)\big] \mathbb{E}\big[V^{2}_{0}(s',
t')\big]\le c.
\end{split}
\end{equation}
Similar to the argument in (\ref{Eq:2S5}) and (\ref{Eq:2S6}), we derive 
from  (\ref{Eq:2P2}) that the integral $\mathcal{K}$ in (\ref{Eq:2P1}) is,
up to a constant,  bounded from above by 
\begin{equation}\label{Eq:2P3}
\arraycolsep=1pt\begin{array}{ll} 
&\displaystyle  \int_{(I \times J)^{2}}\bigg[\sum_{j\in S}\Big(\min\big\{|t_j-t'_j|^{H_j}, t_j^{H_j}\big\}+
\min\big\{|s_j-s'_j|^{H_j}, s_j^{H_j}\big\}\Big)\\
&\displaystyle + \sum_{j\in S^c}\min\big\{|s_j-s_j|^{H_j},
|s_j-s'_j|^{H_j}, |s_j-t'_j|^{H_j}, s_j^{H_j}\big\}\bigg]^{-(d+2) } dsdtds'dt'.
\end{array}
\end{equation}
Again exactly like what we did in the proof of (\ref{Eq:prn2}), we
can show that the integral in (\ref{Eq:2P3}) is finite provided
$2\sum_{j\in S}1/H_j+\sum_{j\in S^c}1/H_j>d+2$. This  proves Part
(iii).

Since the function $\mathbb{E}(V_{0}(s, t) V_{0}(s', t'))$ is uniform continuous
for $(s,\, t, \, s',\,t') \in (I \times J)^2 $ and (\ref{Eq:2S5}) holds, there
exist positive constants $\delta$ and $c$ such that
$\mathbb{E}\big(V_{0}(s, t) V_{0}(s', t')\big) \ge c$
for all  $(s,\,t, \, s',\,t') \in (I \times J)^2 $ such that
$|(s, t) - (s', t')| \le \delta.$
Hence the proof of Part (iv) is quite similar to the proof of Part (ii).
We leave the details to the interested reader. \hfill{$\square$}

\vskip 0.2cm

Part (ii) and (iv) in Theorem \ref{Thm:SmoothSelfO2} can be improved if
we have more information on the dependence structure of $V_0 (s, t) =
X_0(s) - X_0(t)$, as shown by the following theorem.

\begin{thm}\label{Thm:SmoothSelfO2b}
If, in addition to the assumptions of Theorem \ref{Thm:SmoothSelfO2},
$X_0$ satisfies the following condition:
\begin{itemize}
\item[{\rm (C4)}]\ There exists a positive
constant $c_{15}$ such that for all $(s, t),\, (s', t') \in I \times J$,
\begin{equation}\label{Eq:ConV2}
\begin{split}
&{\rm Var}\big(X_0(s)- X_0(t)|X_0(s')- X_0(t')\big)\\
&\qquad \le c_{15} \bigg( \sum_{j\in S} \big(|t_j-t'_j|^{2H_j} +
|s_j-s'_j|^{2H_j}\big)
 + \sum_{j\in S^c}|s_j-t_j|^{2H_j}\bigg).
\end{split}
\end{equation}
\end{itemize}
Then the following statements hold:
\begin{description}
\item {\rm (i)}\  $L_S\left(X^H,I\times J\right)\in L^2(\Omega, \mathbb{P})$ if and only if
$2\sum_{j\in S}\frac1{H_j}+\sum_{j\in S^c}\frac1{H_j}>d.$
\item {\rm (ii)}\ $L_S\left(X^H,I\times J\right)\in D_1$ if and only if
$2\sum_{j\in S}\frac1{H_j}+\sum_{j\in S^c}\frac1{H_j}>d+2.$
\end{description}
\end{thm}

\begin{rem}
Observe that Condition (C4) is automatically satisfied if $S = \emptyset$.
If $X_0 = \{X_0(t), \,t \in \R^N\}$ is an ``additive
fractional Brownian motion'' defined by
\[
X_0(t) = B^{H_1}(t_1) + \cdots + B^{H_N}(t_N), \quad \forall\, t \in \R^N,
\]
where $B^{H_1}, \ldots,  B^{H_N}$ are independent fractional Brownian
motions with indices $H_1, \ldots, H_N$, respectively. Then it is easy to
see that Condition (C4) is satisfied. When $X_0$ is the Brownian sheet,
then by using the independence of increments over intervals, one can
check that (C4) also holds.
\end{rem}

\noindent{\it Proof of Theorem \ref{Thm:SmoothSelfO2b}}
Sufficiencies of the condition in (i) and (ii) have been proved in
Theorem \ref{Thm:SmoothSelfO2}. Note that (C4)
and (\ref{Eq:2S5}) imply
\begin{equation}\label{Eq:2S7b}
  \mathcal{J}\ge  \int_{(I\times J)^{2}}\frac{dsdtds'dt'}
  {\big[\sum\limits_{j\in S}\left(|s_j-s'_j|^{H_j}+|t_j-t'_j|^{H_j}\right)
 + \sum\limits_{j\in S^c}|s_j-t_j|^{H_j}\big]^{d}}.
\end{equation}
By applying Lemma 2.1, it can be verified that the last integral diverges when
$2\sum_{j\in S}\frac1{H_j}+\sum_{j\in S^c}\frac1{H_j}\le d.$ This proves the
necessity in (i). The proof of necessity in (ii) is similar and is omitted.
\hfill{$\square$}

\vskip 0.2cm

Finally, we consider Case (iii), e.g.  the two compact intervals
$I$ and $J$  are not separated in any direction. So $S=\emptyset$.
Compared with Case (ii), we note that, on one hand, (\ref{Eq:2S5})
fails and, on the other hand, Condition (C4) holds
automatically. For concreteness, we assume that $I=J=[0,\,1]^N.$

\begin{thm}\label{Thm:SmoothSelfO3}
Let $X=\{X(t), t\in \mathbb{R}^{N}\}$ be an $(N,
d)$-Gaussian random field as in Theorem \ref{Thm:SmoothSelfO}.
Then the following statements hold:
\begin{description}
\item {\rm (i)}\  $L_S\left(X,I\times I\right)\in L^2(\Omega, \mathbb{P})$
if and only if $\sum_{j=1}^N\frac1{H_j}>d.$
\item {\rm (ii)}\  $L_S\left(X,I\times I\right)\in D_1$ if
$\sum_{j=1}^N\frac1{H_j}> d+2.$
\item {\rm (iii)}\ $L_S\left(X^H,I\times I\right)\notin D_1$ if
\begin{equation}\label{Eq:CWX}
 \sum_{j=1}^{N}\frac1{H_j}\leq \max\Big\{\frac{d+2} 2,\, \frac{2d} 3\Big\}.
\end{equation}
\end{description}
\end{thm}

Before proving this theorem, we compare its conditions with the results in Imkeller
and Weisz (1995, 1999) and Hu (2001). 
\begin{rem}
\begin{itemize}
\item [(a)]\, When $X$ is the $(N, d)$ Brownian sheet, then our existence
condition in (i) coincides with that in Theorem 1 of Imkeller and Weisz (1995, 1999).
When $X$ is a fractional Brownian motion $B^H = \{B^H(t), t \in \R\}$,  our
condition in (ii) becomes $H (d+2) < 1$, which is stronger than that in
Hu (2001, Theorem 3.2). 
\item [(b)]\, Little has been known about optimal necessary condition
for $L_S\left(X,I\times I\right)\in D_1$ for a Gaussian random field $X$.
Our condition (\ref{Eq:CWX}) is the first general result of this kind.
When $X$ is a fractional Brownian motion $B^H$, (\ref{Eq:CWX}) becomes 
$H \ge \min\{\frac 3 {2d},\, \frac 2 {d+2}\}$, which is the complement of the 
sufficient condition of Hu (2001, Theorem 3.2).
Hence, we have proven that, for fractional Brownian motion $ B^H =
\{B^H(t), t \in \R\}$ in $\R^d$ and $I = [0, 1]$, $L_S\left(B^H ,I\times I\right)\in D_1$ 
if and only if $H < \min\{\frac 3 {2d},\, \frac 2 {d+2}\}$.
\end{itemize}
\end{rem}

\noindent{\it Proof of Theorem \ref{Thm:SmoothSelfO3}}
We prove Parts (i) at first.
Notice that by Lemma \ref{Lem:EY2},  we only need to prove that
\begin{equation}\label{Eq:3S1}
\arraycolsep=1pt\begin{array}{ll} \mathcal{J}=\displaystyle\int_{
I^{4}}\frac{dsdtds'dt'}{[\det{\rm Cov}(V_{0}(s, t), V_{0}(s',
t'))]^{d/2} }<\infty
\end{array}
\end{equation}
if and only if $\sum_{j=1}^{N}\frac1{H_j}>d$.
For any $(s,t), (s',t')\in I\times I$, we use  Condition (C1) to
obtain that
\begin{equation}\label{Eq:3S2}
\begin{split}
\det{\rm Cov}(V_{0}(s, t), V_{0}(s', t'))
&\le \E(V^2_0(s,t))\E(V^2_0(s',t'))\\
&\le c\, \bigg(\sum_{j=1}^{N}|s_j-t_j|^{2H_j}\bigg)\bigg(\sum_{j=1}^{N}|s'_j-t'_j|^{2H_j}\bigg).
\end{split}
\end{equation}
This, together with (\ref{Eq:3S1}), implies
\begin{equation}\label{Eq:3S3}
\arraycolsep=1pt\begin{array}{ll} \displaystyle
\mathcal{J}&\displaystyle\geq \displaystyle c\,
\int_{I^{2}}\frac{dsdt} {\big(\sum_{j=1}^{N}|s_j-t_j|^{2H_j}\big)^{d/2}}
\int_{I^{2}} \frac{ds'dt'} { \big(\sum_{j=1}^{N}|s'_j-t'_j|^{2H_j}\big)^{d/2}}\\
&\displaystyle\geq  c\,
\bigg(\int_{I} \frac{dt} {\big(\sum_{\ell=1}^{N}t_{\ell}^{H_{\ell}}\big)^{d}}\bigg)^{2}.
\end{array}
\end{equation}
By using  Lemma \ref{lem:lemma1}, it is elementary to verify that
the last integral in (\ref{Eq:3S3})   is infinite provided
$\sum_{j=1}^{N}1/H_j \leq d$. Hence, we prove the necessity
of Part (i).

To prove the sufficiency in Parts (i), we apply Condition (C3)
to see that for any $(s,t),\,(s',t')\in I\times I$,
\begin{equation*}\label{Eq:3S4}
\begin{split}
{\rm Var}\left(V_0(s,t)|V_0(s',t')\right)&\ge
{\rm Var}\left(X_0(t)|X_0(s),\,X_0(s'),\,X_0(t')\right)\\
&\ge c\,\sum_{j=1}^{N}\min\big\{|t_j-s_j|^{2H_j},
|t_j-s'_j|^{2H_j}, |t_j-t'_j|^{2H_j}, t_j^{2H_j}\big\}.
\end{split}
\end{equation*}
Moreover, we also have
\begin{equation*}\label{Eq:3S5}
\begin{split}
{\rm Var}\left(V_0(s',t')\right)\geq  {\rm Var}\left(X_0(t')|X_0(s')\right)
\ge c\, \sum_{j=1}^{N}\min\big\{ |t'_j-s'_j|^{2H_j},
{t'}_j^{2H_j}\big\}.
\end{split}
\end{equation*}
Combining the above two inequalities with (\ref{Eq:3S1}) yields
\begin{equation}\label{Eq:3S6}
\arraycolsep=1pt\begin{array}{ll} \mathcal{J}\le &\displaystyle
c\, \int_{I ^{4}}\bigg[\sum_{j=1}^{N}\min\big\{|t_j-s_j|^{H_j},
|t_j-s'_j|^{H_j}, |t_j-t'_j|^{H_j},
t_j^{H_j}\big\}\bigg]^{-d}\\
&\displaystyle \qquad  \qquad\times\bigg[\sum_{j=1}^{N}\min\big\{
|t'_j-s'_j|^{H_j}, {t'}_j^{H_j}\big\}\bigg]^{-d}dsdtds'dt'.
\end{array}
\end{equation}
Similarly to the proof of (\ref{Eq:prn2}), we
integrate  $dt_{1}, \ldots, dt_{N}, dt'_{1},
\ldots, dt'_{N}$ to show that the integral in  (\ref{Eq:3S6}) is
finite provided $\sum_{\ell=1}^{N}1/H_{\ell} > d$. This  proves the
sufficiency of Part (i).

In order to prove Part (ii), by Lemma \ref{Lem:CY2}, it suffices to verify
that if $\sum_{j=1}^N\frac1{H_j}>d+2$, then
\begin{equation}\label{Eq:3P1}
\mathcal{K}=\displaystyle\int_{I^4}\frac{\big[\E\left(V_{0}(s, t)
V_{0}(s', t')\right)\big]^2}{[\det{\rm Cov}(V_{0}(s, t), V_{0}(s',
t'))]^{\frac{d+2}{2}} }\,ds\,dt\,ds'\, dt'<\infty.
\end{equation}

For any $(s,t), (s',t')\in I\times I$, the Cauchy-Schwarz
inequality and (C1) imply
\begin{equation}\label{Eq:3P2}
\begin{split}
&\big[\mathbb{E}(V_{0}(s, t) V_{0}(s', t'))\big]^{2}\le
\mathbb{E}\big[V^{2}_{0}(s, t)\big] \mathbb{E}\big[V^{2}_{0}(s',
t')\big]\le c.
\end{split}
\end{equation}
Similar to the argument in Part (i), we derive from  (\ref{Eq:3P2}) that
\begin{equation}\label{Eq:3P3}
\arraycolsep=1pt\begin{array}{ll} \mathcal{J}\le &\displaystyle
c\, \int_{I^{4}}\bigg[\sum_{j=1}^{N}\min\big\{|t_j-s_j|^{H_j},
|t_j-s'_j|^{H_j}, |t_j-t'_j|^{H_j},
t_j^{H_j}\big\}\bigg]^{-(d+2)}\\
&\displaystyle \qquad \qquad \times\bigg[\sum_{j=1}^{N}\min\big\{
|t'_j-s'_j|^{H_j}, {t'}_j^{H_j}\big\}\bigg]^{-(d+2)}dsdtds'dt'.
\end{array}
\end{equation}
Again  we integrate in the order of $dt_{1}, \ldots, dt_{N}, dt'_{1},
\ldots, dt'_{N}$ to show that the integral in  (\ref{Eq:3P3}) is
finite provided $\sum_{\ell=1}^{N}1/H_{\ell} > d+2$. This proves
(\ref{Eq:3P1}) and hence Part (ii).

Finally we prove Part (iii). By taking two disjoint sub-intervals and
argue as in the proof of Theorem \ref{Thm:SmoothSelfO}, one can see
easily that if $2 \sum_{j=1}^{N}\frac1{H_j}\leq d+2$, then the integral
$\mathcal{K} $ in (\ref{Eq:3P1}) diverges and, consequently,
$L_S\left(X^H,I\times I\right)\notin D_1$.

It remains to show that the integral  $\mathcal{K} $ also diverges if
$3\sum_{j=1}^{N}\frac1{H_j}\leq 2d$. To this end, we write
$\rho(s, t) = \sqrt{\E \big(V_{0}(s, t) ^2\big)}$. It will be useful
to note that $\rho(s, t)$ is a pseudo-metric on $\R^{2N}$. Since
\[
\E \big(V_{0}(s, t) - V_{0}(s', t') \big)^2\big) \le 2 \big(\rho(s, t)^2 + \rho(s', t')^2\big),
\]
we see that, if $\rho(s, s') \le \frac 1 2 \rho(s, t)$ and $\rho(t, t')
\le \frac 1 2 \rho(s, t)$, then
\begin{equation}\label{Eq:CWX2}
\begin{split}
\E\left(V_{0}(s, t)  V_{0}(s', t')\right) &= \frac 1 2
\Big[\E \big(V_{0}(s, t) ^2\big) + \E \big(V_{0}(s', t') ^2\big)\\
&\qquad \qquad  - \E \big((V_{0}(s, t) - V_{0}(s', t'))^2\big)\Big] \\
& \ge \frac 1 2 \E \big(V_{0}(s', t') ^2\big).
\end{split}
\end{equation}
Let $B_\rho(s, t) = \{(s', t') \in I^2: \rho(s, s') \le \frac 1 2 \rho(s, t),\,
\rho(t, t') \le \frac 1 2 \rho(s, t)\}$.
It follows from  (\ref{Eq:3P1}), (\ref{Eq:3S2}) and (\ref{Eq:CWX2}) that
%
\begin{equation}\label{Eq:3Q1}
\begin{split}
\mathcal{K} &\ge c\, \int_{I^2} \frac{ ds\,dt} {\rho(s, t)^{d+2}}
\int_{B_\rho(s, t)} \frac{ds'\, dt' } { \rho(s', t')^{d-2} }
\ge \int_{I^2} \frac{ds\,dt} {\rho(s, t)^{2(d - Q)} },
\end{split}
\end{equation}
where $Q = \sum_{j=1}^N \frac 1 {H_j}$. In obtaining the last inequality,
we have used the fact that $\rho(s', t') \le 2 \rho(s, t)$ for all $(s', t')
\in B_\rho(s, t)$, and the Lebesgue measure of $B_\rho(s, t)$ is
$c\, \rho(s, t)^{2Q}$. Under Conditions (C1), $\rho(s, t) \le c_1
\sum_{j=1}^N |s_j-t_j|^{H_j}$ for all $s, t \in I^N$.
We can apply Lemma 2.1 to show that the last integral in (\ref{Eq:3Q1})
diverges if and only if $Q \le 2(d-Q)$. This proves $\mathcal{K}= \infty$
when $3 Q \le 2d$. The proof of Theorem \ref{Thm:SmoothSelfO3} is finished.
\hfill{$\square$}

\vskip 0.2cm

The following are concluding remarks. 
\begin{rem}\, 
\begin{itemize}
\item[(a)]  It is known that Conditions (C1) and (C3) are satisfied by
a large class of Gaussian random fields including $N$-parameter fractional Brownian motion [Pitt (1978)],
fractional Brownian sheets [Wu and Xiao (2007), Xiao (2009)] and stochastic heat equation driver by
space-time Gaussian noises [Dalang, et al. (2015), Tudor and Xiao (2015)]. Hence Theorems 
\ref{Thm:SmoothSelfO} and \ref{Thm:SmoothSelfO2} can be applied directly to these
Gaussian random fields. However, despite the conditions given by Theorems 
\ref{Thm:SmoothSelfO2b}  and \ref{Thm:SmoothSelfO3}, 
the problem for finding necessary and sufficient conditions
for $L_S\left(X,I\times I\right)\in D_1$ is still open for a general Gaussian random field.
It would be interesting to solve this problem.
\item[(b)] Another interesting question is to remove the i.i.d. assumption on the coordinate random fields 
$X_1, \ldots X_d$ in (1.1). While the results of this paper can be extended to Gaussian random fields 
with independent, but non-identically distributed components, it seems more difficult to remove the 
independence assumption.  Some preliminary results have been proved by Eddahbi, et al. (2005, 2007)
for vector-valued fractional Brownian sheets, but their conditions may not be optimal.
\end{itemize}
\end{rem}

\noindent\\

\noindent\bf{\footnotesize Acknowledgements}\quad\rm
{\footnotesize Research of Z. Chen and D. Wu was partially supported
by the National Natural Science Foundation of China (Grant No. 11371321).
Research of Y. Xiao was partially
supported by the NSF Grants  DMS-1307470 and DMS-1309856.}\\[4mm]
\space

\vskip 0.1cm

\noindent{\bbb{References}}

\end{document}